\pdfoutput=1
\documentclass[a4paper]{amsart}
\usepackage[english]{babel}

\usepackage{amssymb}
\usepackage{amsmath,amsthm}
\usepackage{mathtools}
\usepackage{dsfont}
\usepackage{nicefrac}
\usepackage{cases}
\usepackage{bm}
\usepackage{enumitem}

\usepackage[T1]{fontenc}
\usepackage{lmodern}
\usepackage{microtype}

\usepackage{empheq}

\usepackage[a4paper, inner=3cm,outer=3cm,top=3.3cm,bottom=3cm]{geometry}

\usepackage{xspace}

\usepackage[autostyle]{csquotes}

\usepackage{tikz}
\usepackage{tikz-cd}
\usetikzlibrary{matrix,arrows,patterns,decorations.pathreplacing,svg.path}
\usepackage{pgfplots}

\usepackage{todonotes}

\usepackage[colorlinks, bookmarksdepth=3, citecolor=blue,linkcolor=blue]{hyperref}

\usepackage{pdfcomment}

\makeatletter
\@ifpackageloaded{hyperref}{%
\DeclareRobustCommand{\SkipTocEntry}[5]{}}{%
\DeclareRobustCommand{\SkipTocEntry}[4]{}}
\makeatother

\newenvironment{pf}{\proof[\proofname]}{\endproof}
\theoremstyle{plain}
\newtheorem{Th}{Theorem}[section]
\newtheorem{Prop}[Th]{Proposition}
\newtheorem{Cor}[Th]{Corollary}
\newtheorem{Lem}[Th]{Lemma}

\numberwithin{equation}{section}
\numberwithin{figure}{section}
\theoremstyle{definition}
\newtheorem{Def}[Th]{Definition}
\newtheorem{Ex}[Th]{Example}
\newtheorem{Rem}[Th]{Remark}

\newcommand{\cone}{{\mathrm{cone}}}
\newcommand{\rk}{{\mathrm{rk}}}

\newcommand{\Prt}{\operatorname{Prt}}

\newcommand{\Km}{{\mathcal{K}}}

\newcommand{\PMV}{{\mathcal{PV}}}
\newcommand{\MV}{{\mathcal{V}}}

\newcommand{\cal}[1]{\mathcal{#1}}

\newcommand{\R}{\mathbb R}

\newcommand{\N}{\mathbb N}

\newcommand{\RP}{\mathbb{RP}}
\newcommand{\G}{\Gamma}

\newcommand{\cA}{\cal A}

\newcommand{\cH}{\cal H}

\newcommand{\cK}{\cal K}
\newcommand{\cL}{\cal L}
\newcommand{\cM}{\cal M}

\newcommand{\cZ}{\cal Z}

\newcommand{\ph}{\varphi}

\newcommand{\la}{\langle}
\newcommand{\ra}{\rangle}

\newcommand{\1}{\mathbf{1}}

\DeclareMathOperator{\V}{V}
\DeclareMathOperator{\PV}{PV}
\DeclareMathOperator{\vv}{v}

\DeclareMathOperator{\Vol}{Vol}

\DeclareMathOperator{\Gr}{Gr}

\newcommand{\Abs}{\operatorname{Abs}}

\newcommand{\spn}{\operatorname{span}}

\newcommand{\im}{\operatorname{Im}}
\newcommand{\Ker}{\operatorname{Ker}}

\newcommand{\multiset}[2]{
\left.\mathchoice
  {\left(\kern-0.48em\binom{#1}{#2}\kern-0.48em\right)}
  {\big(\kern-0.30em\binom{\smash{#1}}{\smash{#2}}\kern-0.30em\big)}
  {\left(\kern-0.30em\binom{\smash{#1}}{\smash{#2}}\kern-0.30em\right)}
  {\left(\kern-0.30em\binom{\smash{#1}}{\smash{#2}}\kern-0.30em\right)}
\right.}

\makeatletter
\renewcommand{\boxed}[2][\fboxsep]{{%
  \setlength{\fboxsep}{#1}\fbox{\m@th$\displaystyle#2$}}}
\makeatother

\newcommand{\rs}[1]{Section~\ref{S:#1}}
\newcommand{\rl}[1]{Lemma~\ref{L:#1}}
\newcommand{\rp}[1]{Proposition~\ref{P:#1}}

\newcommand{\rr}[1]{Remark~\ref{R:#1}}

\newcommand{\re}[1]{(\ref{e:#1})}

\newcommand{\rt}[1] {Theorem~\ref{T:#1}}
\newcommand{\rd}[1]{Definition~\ref{D:#1}}
\newcommand{\rf}[1]{Figure~\ref{F:#1}}

\makeatletter
\@namedef{subjclassname@2020}{%
  \textup{2020} Mathematics Subject Classification}
\makeatother

\title[Mixed volumes of zonoids and the absolute value of the Grassmannian]{
Mixed volumes of zonoids and\\ the absolute value of the Grassmannian}
\author{Gennadiy~Averkov}
\address{Fakult\"at 1, BTU Cottbus-Senftenberg, Platz der Deutschen Einheit 1, 03046 Cottbus, Germany}
\email{averkov@b-tu.de}
\author{Katherina von Dichter}
\email{vondicht@b-tu.de}
\author{Simon~Richard}
\address{Department of Mathematics and Statistics, Cleveland State University,  2121 Euclid Ave, Cleveland, Ohio, 44115 USA}
\email{j.s.richard@vikes.csuohio.edu}
\author{Ivan~Soprunov}
\email{i.soprunov@csuohio.edu}

\begin{document}
\selectlanguage{english}

\date{}


\maketitle

\begin{abstract}
	Zonoids are Hausdorff limits of zonotopes, while zonotopes are convex polytopes defined as the Minkowski sums of finitely many segments. We present a combinatorial framework that links the study of mixed volumes of zonoids (a topic that has applications in algebraic combinatorics) with the study of the absolute value of the Grassmannian, defined as the image of the Grassmannian under the coordinate-wise absolute value map. We use polyhedral computations to derive new families of inequalities for $n$ zonoids in dimension $d$, when $(n,d)=(6,2)$ and $(6,3)$. 
Unlike the classical geometric inequalities, originating from the Brunn-Minkowski and Aleksandrov-Fenchel inequalities, the inequalities we produce have the special feature of being Minkowski linear in each of the $n$ zonoids they involve.
\end{abstract}

\section{Introduction} The mixed volume assigns a non-negative number to any collection of $d$ convex bodies in $\R^d$. This notion goes back to the work of Minkowski and Aleksandrov on Brunn-Minkowski theory of convex bodies at the beginning of 20th century. Many metric functions of a convex body, such as the volume, the surface area, the mean width, and the area of a projection are instances of the mixed volume. 
%Surprisingly, over the past hundred years mixed volumes have made an appearance in almost all areas of pure and applied mathematics: combinatorics, algebraic geometry, optimization, as well as information theory and stochastic geometry. 
%To a large extent, convex geometry is concerned with finding geometric inequalities. 
Thus, a search for inequalities between mixed volumes is an intrinsic problem in convex geometry. The most famous example of such an inequality is  the Aleksandrov-Fenchel inequality. It is a central result in the Brunn-Minkowski theory %. For example, the isoperimetric inequality is an instance of the Aleksandrov-Fenchel inequality, 
which has deep connections to  algebraic geometry via the \emph{Hodge Index Theorem} and to combinatorics, where it explains the log-concavity phenomenon in many problems; see \cite{adiprasito2017hodge,adiprasito2018hodge,branden2020lorentzian,huang2025realizations,khovanskii1978newton,teissier1979theoreme}. Most classical inequalities, e.g. Fenchel's inequalities and Shephard's determinantal inequalities can be derived from the Aleksandrov-Fenchel inequality. Yet, they do not provide a complete system of polynomial inequalities in most cases, where by a complete system we mean a system of inequalities that exhaustively describes the relations between the functionals involved into it. 
 The question of providing a complete system of inequalities for mixed volumes of a given number of bodies $n$ in a given dimension $d$ goes back to the work of Heine in 1938 \cite{Heine38} where he gave such a system for $n=3$ bodies in dimension $d=2$. The inequalities he produced were later extended by Shephard \cite{Shephard60} to the case of arbitrary dimension. However, Shephard's inequalities give a complete system only for $n=2$ and are not sufficient for $n>2$. All the other cases of Heine-Shephard problem remain currently open. 

In \cite{AS23} Averkov and Soprunov addressed the smallest unknown case of $n=4$ bodies in dimension $d=2$. They introduced Plücker-type inequalities for the six pair-wise mixed areas of the four bodies and proved that they form a complete system for these six functionals. Note that, unlike in the original Heine-Shephard problem, the four areas are discarded and only the ``purely mixed'' areas are considered. These inequalities are closely related to the Plücker relation for the Grassmannian of 2-planes in 4-space. In this paper we explore this connection further for a specific class of convex bodies, the zonoids.

Recall that zonotopes are Minkowski sums of segments and zonoids are their limits in the Hausdorff metric. 
 Zonotopes are naturally related to matroids. For example, in \cite{Stan} Stanley used the Aleksandrov-Fenchel inequality for the mixed volume of zonotopes to show general log-concavity behavior of certain sequences that count the number of bases of a regular matroid. More recently, mixed volumes of zonoids and zonotopes have been studied in analysis \cite{FZ23} and probability (random determinants) \cite{ZonoidAlgebra}. Thus, a systematic study of inequalities between mixed volumes of zonoids is of particular interest across different research areas. 

\subsection{Summary of the results} Given $n,d\in\N$ with $n$ divisible by $d$, we are looking for homogeneous polynomial inequalities of degree $n/d$ relating ${n\choose d}$ values of  mixed volumes $\V(Z_{i_1},\dots, Z_{i_d})$, for $1 \le i_1 < \cdots < i_d \le n$, which are valid for all collections $Z=(Z_1,\dots, Z_n)$ of $n$ zonoids 
in $\R^d$. If we think of these mixed volumes as variables (which can be evaluated at a collection $Z$)
then every such inequality can be viewed as a linear inequality in its monomials in these variables. 
Geometrically, this means that we have a ``monomial'' map $\Phi$ that sends every collection 
$Z$  of $n$ zonoids to the vector of monomials evaluated at $Z$, and we are interested in the conic hull $C_{n,d}$ of the image of this map $\Phi$. The facets of $C_{n,d}$ give us, in some sense, the ``best'' inequalities of the multilinear type.  To digest this somewhat technical construction one may consider, for example, the case $(n,d)=(4,2)$, where the corresponding cone 
\[
		C_{4,2} = \cone \{ ( \V_{12} \V_{34} , \V_{13} \V_{24} , \V_{14} \V_{23} ) \colon \V_{ij} = \V(Z_i,Z_j), \ Z_1,\ldots, Z_4  \ \text{zonoids in } \ \R^2\} \subset \R^3.
\]
The ambient space of $C_{4,2}$ has dimension $3$, which is the number of $2+2$ partitions of $[4]$. Linear inequalities valid on $C_{4,2}$ are in one-to-one correspondence with the multilinear inequalities for four zonoids $Z_1,\ldots,Z_4$ in the plane. 

The distinctive feature of the inequalities we consider is their linearity in each $Z_i$ with respect to Minkowski addition, see \rs{multilinear}. This implies that the verification of such inequalities is equivalent to verifying their validity for line segments, as every zonotope is generated by line segments and every zonoid is a limit of zonotopes (whereas mixed volumes are continuous functionals). This observation has two consequences. First, it connects our problem to combinatorics of point configurations in real projective spaces, see \rs{proj-view}. It turns out that there are certain ``rigid'' point configurations which (in the cases we consider) describe the generators of the cone
$C_{n,d}$ in the cases $d=2$ and $(n,d) = (6,3)$. We prove this in \rs{rigid}. How to describe the generators of $C_{n,d}$ in general remains an open problem. 

Second, this allows deriving such inequalities  via polyhedral computations with the help of a computer. In Sections \ref{S:framework-4-2}-\ref{S:framework-6-3} we show that when
$(n,d)=(4,2), (6,2)$ and $(6,3)$ the conic hull $C_{n,d}$ is a polyhedral cone of dimension 3, 15, and 10,
with 3, 25, and 30 generators, respectively. In the smallest case $(n,d)=(4,2)$ the cone has three facets and
we recover the three Plücker-type inequalities introduced in \cite{AS23}. Expectedly, the number of facets in the other two cases is much larger: 975 for $(n,d)=(6,2)$ and 130 for $(n,d)=(6,3)$. However, up to the action of the symmetric group permuting the six zonoids, there are only eight new inequalities in the former, and two in the latter case.
Our computations were done using {\it SageMath} \cite{sagemath}. The code is available at \url{https://github.com/jsimonrichard/mixed-volumes-n-6-d-2-3}.

There are several phenomena we have observed by studying these more complicated cases, beyond the case $(n,d)=(4,2)$.

\begin{enumerate}
\item In contrast to the case $(n,d)=(4,2)$, the ``best'' Minkowski multilinear inequalities no longer give a complete system for the ``purely mixed'' volumes when the values of $n$ and $d$ are higher, see \rp{2-cells-(2,2,2)}.

\item Note that in dimension $d=2$, zonoids are just centrally-symmetric convex bodies. It was shown in \cite{AS23} that for $(n,d)=(4,2)$ there is no difference between complete systems in
the centrally-symmetric and in the general case. However, for larger values of $n$ and $d$ such a difference does emerge. For example, for $(n,d)=(6,2)$ we find inequalities which hold in the centrally-symmetric case, but not in the general case, see \rr{discrepancy}. It would be interesting to fully explore this discrepancy in this case and for larger values of $n$ and $d$. 

\item In our studies of $C_{n,d}$, by letting $n$ and $d$ grow, we have observed a complexity blow up, at an early stage, along with a qualitative structural change from polyhedrality to general convexity later on. When $n$ or $d$ are increased, the number of defining inequalities of $C_{n,d}$ grows rapidly and, at a certain threshold, the cone $C_{n,d}$ passes from the state of being a polyhedron to the state of being a non-polyhedral closed semialgebraic set. In \rs{nonpolyhedrality} we show that $C_{n,2}$ is not polyhedral for every even $n \ge 8$. Our studies suggest that finding complete systems of inequalities without computer assistance becomes infeasible for most choices of $(n,d)$.
This underscores the importance of computer-assisted methods in our context.
Traditionally, geometric inequalities have been discovered and verified ``by hand.'' To the best of our knowledge, there is no prior work in convex geometry in which geometric inequalities were determined using computer-assisted tools.

\end{enumerate}

\subsection{Connection to absolute value of the Grassmannian} On the algebraic side, our problem is closely related to studying the {\it absolute value} of the real Grassmannian  $\Gr_{d,n}$ which is the set %\subset \R^{{n\choose d}}$ %of $d$-planes in $\R^n$. 
\[
	\Abs( \Gr_{d,n} )  = \{ \bigl( | \det(u_{i_1},\ldots,u_{i_d}) | \bigr)_{1 \le i_1 \le \cdots \le i_d \le n} \colon (u_1,\ldots,u_n) \in \R^d \}\subset \R_{\geq 0}^{{n\choose d}}.
\]
%Here $\Abs:\R^{{n\choose d}}\to \R_{\geq 0}^{{n\choose d}}$ denotes the component-wise absolute value map. 
In fact, the inequalities for mixed volumes of zonotopes produced using our method are the facet inequalities for the conic hull of $\Abs(\Gr_{d,n})$, see \rs{abs-grass} for details.

Points of $\Abs(\Gr_{d,n})$ store some metric information about a vector configuration $(u_1,\ldots,u_d)$ (the absolute values of the determinants), but do not keep track of the orientation of the configuration.
Such a setting looks very similar to a first step in deriving a tropicalization, where in order to tropicalize an algebraic variety $V \subseteq \mathbb{F}^n$ over a field $\mathbb{F}$ with a valuation, one applies the valuation to the points of $V$ in the component-wise fashion and then takes the closure of the resulting set. The tropical Grassmannian was studied in \cite{speyer2004tropical}. Whether there exist any concrete connections between the absolute and the tropical Grassmannian is still to be explored. It would also be interesting to study potential connections of $\Abs(\Gr_{d,n})$ to the positive Grassmannian $\Gr_{d,n}^{\ge 0} := \Gr_{d,n} \cap \R_{\ge 0}^{\binom{n}{d}}$, see \cite{postnikov2018positive,williams2021positive}. The connection that is immediately clear is that  $\Gr_{d,n}^{\ge 0}$ is a subset of $\Abs(\Gr_{d,n})$, so that any inequality valid on $\Abs(\Gr_{d,n})$ is also valid on $\Gr_{d,n}^{\ge 0}$.

\subsection*{Acknowledgments} Ivan Soprunov is grateful to Darren Gerrity and Kyle Klingensmith for helpful discussion regarding the result of \rt{nonpolyhedral}. The authors thank the anonymous referees for valuable suggestions. Ivan Soprunov is supported by the AMS-Simons Travel Grant. Katherina von Dichter is supported by the Postdoc Network Brandenburg (Germany) through the project \emph{Geometrische Ungleichungen für gemischte Volumina} (English: Geometric Inequalities for Mixed Volumes).

\section{Preliminaries}
For a positive integer $m$, we use $[m]$ to denote $\{1,\ldots,m\}$. 
By $e_1,\ldots,e_d$ we denote the standard basis vectors of $\R^d$. 
%\subsection{Mixed volumes}
Let $\Km_d$ denote the set of all convex bodies, that is, non-empty compact convex subsets of $\R^d$.
For $K\in\Km_d$ let $\Vol_d(K)$  denote the Euclidean $d$-dimensional volume of $K$. Given two subsets $A,B$ of $\R^d$, their \emph{Minkowski sum} $A+B$ is the vector sum $A+B=\{a+b : a\in A, b\in B\}$. 

\subsection{Mixed volume}
The {\it mixed volume} $\V(K_1,\dots, K_d)$ is the unique symmetric and multilinear (with respect to Minkowski addition) function of $d$ bodies $K_1,\dots, K_d$ in $\Km_d$, which coincides with $\Vol_d$ on the diagonal, i.e., 
$\V(K,\dots, K)=\Vol_d(K)$ for any $K\in\Km_d$. One can write an explicit ``polarization'' formula for the mixed volume \cite[Section 5.1]{Schneider2014}:
\begin{equation}\label{e:polarization}
\V(K_1,\dots,K_d)=\frac{1}{d!}\sum_{\ell=1}^d(-1)^{d+\ell}\!\sum_{1 \le i_1<\dots<i_\ell \le d}\Vol_d(K_{i_1}+\dots+K_{i_\ell}).
\end{equation}
Basic properties of the mixed volume include non-negativity, invariance under independent translations of the $K_i$, as well as monotonicity with respect to inclusion of the bodies. 

A very special case which links
mixed volumes and matroids is when the bodies are line segments. In this case the mixed volume
equals the volume of the parallelotope spanned by the segments, divided by $d!$, by \re{polarization}. More
explicitly, let $L_1,\dots, L_d$ be line segments in $\R^d$, where $L_i=[0,u_i]$ for some $u_i\in\R^d$. Then we have 
\begin{equation}\label{e:det}
\V(L_1,\dots, L_d)=\frac{1}{d!}\left|\det(u_1,\dots,u_d)\right|.
\end{equation}

A {\it zonotope} is a convex polytope which is the Minkowski sum of line segments. A convex body which is a limit of zonotopes in Hausdorff metric is called a {\it zonoid}, see \cite[Section 3.5]{Schneider2014}. We use  $\cZ_d$ to denote the space of all zonoids in $\R^d$. 
It is well-known that every centrally-symmetric convex polygon in $\R^2$ is the Minkowski sum of finitely many segments, i.e. is a zonotope. Therefore, in dimension $2$ the space of zonoids $\cZ_2$ coincides with the space of all centrally-symmetric convex bodies.

\subsection{Projective viewpoint}\label{S:proj-view} 
In this paper we study homogeneous polynomial inequalities 
involving mixed volumes of $n$ zonotopes in $\R^d$. Our inequalities are Minkowski linear in each zonotope, which implies that their verification boils down to verifying them on 
$n$ line segments in $\R^d$. By translation invariance of the mixed volume, we can 
identify an $n$-tuple of line segments $L=(L_1,\dots, L_n)$ with an $n$-tuple of non-zero vectors $U=(u_1,\dots, u_n)$ in $\R^d$ all lying in a half-space, which can be chosen to be the half-space of vectors with non-negative last coordinate. This allows us to associate $U=(u_1,\dots, u_n)$ with an $n$-tuple of points $A=(a_1,\dots, a_n)$
in the real projective space $\RP^{d-1}$, where the points may coincide if the corresponding vectors are proportional. Finally, if two $n$-tuples $U$ and $U'$ coincide up to a non-singular linear transformation
then  the corresponding $n$-tuple of points $A$ in $A'$ are projectively equivalent. This amounts to 
rescaling the mixed volumes by a positive scalar which would preserve our homogeneous inequalities.
This projective viewpoint is helpful when deriving the inequalities. Often, when studying geometry of
$n$-tuples $A$ in $\RP^{d-1}$ we will forget the order and view $A$ as a {\it configuration} (a multiset)
of points in $\RP^{d-1}$. We denote a configuration by listing the elements with multiplicities. For example, 
$\{2\cdot a,\, 3\cdot b,\, c\}$ denotes a configuration with a double point $a$, a triple point $b$, and a simple point $c$. We also abbreviate $2\cdot\{a,b\}$ for a configuration with two double points.

\subsection{Mixed volume configuration spaces}
Consider an $n$-tuple $K=(K_1,\dots, K_n)$ of bodies in $\Km_d$. We use 
${[n]\choose d}$ to denote the set of all subsets of $[n]$ of size $d$ and 
$\multiset{[n]}{d}$ to denote the set of all multisets on $[n]$ of size $d$. 
%Each $I=\{m_1\cdot 1, \dots, m_n\cdot n\}\in \multiset{[n]}{d}$ 
Each $I\in \multiset{[n]}{d}$, where $I = \{ \underbrace{1,\ldots,1}_{m_1},\ldots, \underbrace{n,\ldots,n}_{m_n}\}$,
defines the value of the mixed volume 
$$\V_I(K):=\V( \underbrace{K_1,\ldots,K_1}_{m_1},\ldots,\underbrace{K_n,\ldots,K_n}_{m_n}).$$
Then the {\it mixed volume configuration} corresponding to $K$ is the vector of mixed volumes over all multisets $\V(K):=\left(\V_I(K) : I\in\multiset{[n]}{d}\right)\in\R^{d+n-1\choose d} $.

\begin{Def}\label{D:MV-space} 
Let $\cA\subset\cK_d$ be a set of convex bodies in $\R^d$.
Define the {\it mixed volume configuration space} $\MV(n,\cA)$ to be the set of 
all mixed volume configurations $\V(K)$ over all $n$-tuples $K\in \cA^n$.
\end{Def}

\begin{Ex}\label{Ex:22} For $n=d=2$ the mixed volume configuration corresponding to $K=(K_1,K_2)$ is the
vector $\V(K)=\left(\V(K_1,K_1),\V(K_1,K_2),\V(K_2,K_2)\right)$. Note that the first and the last entries are
the areas $\Vol_2(K_1)$ and $\Vol_2(K_2)$, respectively. According to the Minkowski inequality, 
$$\V(K_1,K_2)^2\geq \Vol_2(K_1)\Vol_2(K_2).$$
This implies that $\MV(2,\cK_2)\subseteq \{(x,y,z)\in\R^3 : y^2\geq xz,\, x\geq 0,\, y\geq 0,\, z\geq 0\}$.
It is not hard to see that the above inclusion is, in fact, equality. 
%In particularly, $\MV(2,\cK_2)$ is a {\it semialgebraic set}, that is, a set which can be described by a boolean combination of polynomial inequalities. 
\end{Ex}

For the purpose of this paper we will assume that $n\geq d$ and restrict ourselves to subsets $I\subset[n]$
of size $d$, rather than multisets. Namely, given an $n$-tuple $K\in\cK_d^n$, 
define the {\it pure mixed volume configuration} corresponding to $K$ as the vector of mixed volumes over all subsets $\PV(K):=\left(\V_I(K) : I\in{[n]\choose d}\right)\in\R^{n\choose d} $.

\begin{Def}\label{D:PMV-space} 
Let $\cA\subset\cK_d$ be a set of convex bodies in $\R^d$.
Define the {\it pure mixed volume configuration space} $\PMV(n,\cA)$ to be the set of 
all pure mixed volume configurations $\PV(K)$ over all $n$-tuples $K\in \cA^n$.
\end{Def}

Note that $\PMV(n,\cA)$ is the image of $\MV(n,\cA)$ under projection onto a coordinate subspace.
One readily sees that $\PMV(d,\cK_d)=\R_{\geq 0}$, since a single mixed volume can take an arbitrary non-negative value. Also, one can show that $\PMV(d+1,\cK_d)=\R^{d+1}_{\geq 0}$. The first non-trivial example
occurs for $n=4$ and $d=2$. It was shown in \cite{AS23} that $\PMV(4,\cK_2)$ is a semialgebraic subset of $\R^{6}_{\geq 0}$ given by three quadratic inequalities, called Plücker-type inequalities (see also \rs{framework-4-2} below). The main goal of the paper is to produce similar Minkowski multilinear homogeneous inequalities which hold for $\PMV(n,\cZ_d)$ for larger values of $n$ and $d$.

\subsection{Minkowski multilinear inequalities and the cone $C_{n,d}$}\label{S:multilinear}
Let $K=(K_1,\dots, K_n)\in \cK_d^n$ be an $n$-tuple of convex bodies and, for each $I=\{i_1,\dots, i_d\}\subset [n]$, let $\V_I(K)=\V(K_{i_1},\dots, K_{i_d})$ be the corresponding
pure mixed volume. 
We assume that $n=kd$ for some integer $k$. Then every partition  $I_1|\cdots| I_k=[n]$ into $k$ subsets of size $d$ gives rise to 
a degree $k$ monomial $\V_{I_1}(K)\cdots\V_{I_k}(K)$. We are looking for inequalities of the form 
\begin{equation}\label{e:framework}
\sum_{I_1|\cdots|I_k=[n]}\lambda_{I_1|\cdots|I_k}\V_{I_1}(K)\cdots\V_{I_k}(K)\geq 0,
\end{equation}
where $\lambda_{I_1|\cdots|I_k}\in\R$ and the sum is over all partitions of $[n]$ into $k$ subsets of size $d$. 

To simplify notation, we will often use $\Prt(n,d)$ to denote the set of all partitions
of $[n]$ into $d$-element subsets. Furthermore, we will abbreviate 
\begin{equation}\label{e:abbreviate}
\V_P(K):=\V_{I_1}(K)\cdots\V_{I_k}(K),\quad\text{for }\ P=I_1|\cdots|I_k\in\Prt(n,d)
\end{equation}
when it is convenient. In this notation \re{framework} becomes
\begin{equation}\label{e:framework-short}
\sum_{P\in\Prt(n,d)}\lambda_{P}\V_{P}(K)\geq 0.
\end{equation}

Let $\left( x_I : I\in {[n] \choose d}\right)$ denote the system of coordinates for $\R^{n\choose d}$. Consider the monomial map 
\begin{equation}\label{e:monomial-map}
\ph:\R^{n\choose d}\to\R^{N},\quad \ph(x)= \left(x_{I_1}\cdots x_{I_k}\,:\,  I_1|\cdots|I_k=[n]\right),
\end{equation}
where $N=\frac{n!}{k!d!^k}$ is the number of partitions. We call $\ph$ the {\it partition monomial map}.
Precomposing $\ph$ with the mixed volume functional we obtain the map
sending an $n$-tuple $K=(K_1,\dots, K_n)$ of convex bodies to the vector of all monomials appearing in \re{framework}:
\begin{equation}\label{e:our-map}
\Phi:\cK_d^n\to\R^{N},\quad \Phi(K)= \left(\V_{I_1}(K)\cdots\V_{I_k}(K)\,:\,  I_1|\cdots|I_k=[n]\right).
\end{equation}
%where $N=\frac{n!}{k!d!^k}$ is the number of partitions. 
This map is Minkowski multilinear, that is
\begin{equation}\label{e:linearity}
\Phi(K_1,\dots, \lambda K_i+\mu L_i,\dots,K_n)=\lambda\,\Phi(K_1,\dots, K_i,\dots,K_n)+\mu\,\Phi(K_1,\dots, L_i,\dots,K_n)
\end{equation}
for any $i\in[n]$, $\lambda,\mu\in\R_{\geq 0}$, and $K_1,\dots,K_n,L_i\in\cK_d$. This follows from the corresponding property of the mixed volume and
the fact that in every monomial $\V_{I_1}(K)\cdots\V_{I_k}(K)$ there is exactly one mixed volume containing $K_i$ exactly once. 

The cone $C_{n,d}$ is now defined as 
\begin{align*} 
C_{n,d}  & = \cone  \PMV(n, \cZ_d)
		\\ & =  \cone \{ (\V_P(Z) )_{P \in \Prt(n,d) } \colon Z \in \cZ_d^n \} 	
		\\ & =  \cone \{ (\V_{I_1}(Z) \cdots \V_{I_k}(Z) ) \colon I_1 | \cdots | I_k = [n], \  Z \in \cZ_d^n \} 	
\end{align*} 
%\[
%	 
%	 = \cone \{  \}. 
%\]

\subsection{The absolute value of the Grassmannian}\label{S:abs-grass}
%We now turn to the absolute value of the Grassmannian mentioned in the introduction. 
Recall that the Grassmannian $\Gr_{d,n}$ of $d$-planes in $\R^n$ is embedded into $\RP^{{n \choose d}-1}$ via the Plücker coordinates as follows. Given a $d$-plane $W$ in $\R^n$, let $U$ be a $d\times n$ matrix whose 
rows form a basis for $W$.  Then  
$$\Gr_{d,n}\to \RP^{{n \choose d}-1},\quad W\mapsto [\,\det(U_I) : I\in{[n]\choose d}\,] $$
where $U_I$ is the $d\times d$ submatrix of $U$ indexed by $I$. 
Define the {\it absolute value of the Grassmannian} $\Abs(\Gr_{d,n})$ to be the set of orbits of the points 
$(|\det(U_I)| : I\in{[n]\choose d})$ under the diagonal action of the multiplicative group $\R_{>0}$.
Since the mixed volume of segments is the absolute value of the corresponding determinant (see \re{det}), we
can identify $\Abs(\Gr_{d,n})$ as a subspace of $\PMV(n,\cZ_d)$ (modulo the action of  $\R_{>0}$). In this way, 
we have a description of the cone $C_{n,d}$ as
\[
	C_{n,d} = \cone \, \ph(\Abs(\Gr_{d,n})).  
\]
This description also shows that the set $C_{n,d}$ is a semialgebraic convex cone, in general. This follows from 
the theory of semialgebraic sets (see \cite{bochnak2013real,scheiderer2024course} for more background), as $\Gr_{d,n}$ is a real algebraic variety, $\Abs$ is a semialgebraic map (in fact, a piecewise linear map), $\ph$ is a polynomial map and the conic hull of a semialgebraic set is semialgebraic. 
In particular, every Minkowski multilinear inequality \re{framework} holds for the %pure mixed volume configuration 
space $\PMV(n,\cZ_d)$ if and only if it holds for %the absolute value of the Grassmannian 
$\Abs(\Gr_{d,n})$.

%%%%%%%%%%%%%%%%%%%%%%

\section{Fiber-invariance of the partition monomial map for $d=2$}\label{S:fiber-invariance}

Recall the definition of the pure mixed volume configuration space $\PMV(n,\cA)$ from \rd{PMV-space}.
In the notation introduced in \rs{multilinear}, we have $\ph(\PMV(n,\cA))=\Phi(\cA)$.
Assume $\cA\subseteq\cK_2$ is closed under taking positive linear combinations of the bodies.
The following theorem shows that, restricting to the positive orthant, 
the problems of providing a description for $\PMV(n , \cA)$ and for $\Phi(\cA)$ are equivalent. 

\begin{Th}\label{T:fiber-invariance} Let $\cA\subseteq\cK_2$ be a set of convex bodies closed under taking positive linear combinations. Then for any $x\in\R_{>0}^{n\choose 2}$ we have $x\in\PMV(n , \cA)$ if and only if
$\ph(x)\in\ph(\PMV(n,\cA))=\Phi(\cA)$.
\end{Th}

We will deduce \rt{fiber-invariance} from a general result about fiber-invariance of certain sets under the partition monomial map $\ph$, see \rt{main-fiber}.

In the case of $d=2$ and $n$ even, it is convenient to identify the set
of partitions of $[n]$ into 2-element subsets with the set $\cM_n$ of perfect matchings in the complete graph on $n$ nodes.
We also identify  $\binom{[n]}{2}$ with the set of edges $E_n$ of the complete graph and write $x_e$ for the variable 
corresponding to $e\in E_n$. Then the partition monomial map \re{monomial-map}, restricted to the positive orthant, becomes

\begin{equation}\label{e:ph-map}
\ph:\R_{>0}^{E_n}\to\R_{>0}^{\cM_n},\quad \ph(x)= \left(x_M\,:\,  M\in\cM_n\right),\quad\text{where }\ x_M=\prod_{e\in M}x_e.
\end{equation}

We proceed with a general definition.

\begin{Def} Let $f : X \to Y$ be a map between two sets. A subset $D \subseteq X$ is called {\it fiber-invariant for $f$} if for every $y\in f(D)$ the fiber
$f^{-1}(y) := \{ x \in X \colon f(x) = y\}$ is contained in~$D$.
\end{Def} 

Note that fiber invariance of $D$ with respect to $f$ means that $D$ is the preimage of some set with respect to $D$ or, in other words, one has $D = f^{-1}( f(D))$. Somewhat informally, this means that $f(D)$ keeps enough information for the description of $D$. 

\begin{Ex} 
Consider the set $D = \{ (x,y,z)\in\R^3 : z \ge x^2 + y^2 \}$. Then $D$ is fiber-invariant for the map
$f : \R^3 \to \R^2$ given by $f(x,y,z) =  (x^2 + y^2, z)$. However, $D$ is not fiber-invariant for $g(x,y,z) = (x+y,z)$. To see the latter, observe that the fiber of $(1,1)$ contains the point $(1,0,1)$ which belongs to $D$, but also the point $(-1,2,1)$ which does not belong to~$D$. 
\end{Ex} 

We will need the following simple lemma.

\begin{Lem}\label{L:exact-seq} 
Consider a sequence of vector spaces 
$\begin{tikzcd}   U \arrow[r, "g"] & V \arrow[r, "f"] & W \end{tikzcd}$ which is exact in $V$. Then 
a subset $D\subset V$ is fiber-invariant for $f$ if and only if $D=D+\Ker f=D+\im g$.
\end{Lem}

\begin{pf} By exactness in $V$ we have $\Ker f=\im g$. Now the lemma
follows immediately from the definition of fiber-invariance and the fact that for any $x\in V$
the fiber of $y=f(x)$ is $f^{-1}(y)=x+\Ker f$.
\end{pf}

It follows from the definition that if $D$ is fiber-invariant for $f:X\to Y$, then checking whether $x \in X$ belongs to $D$ is equivalent to checking whether $f(x)$ belongs to $f(D)$. 
Indeed, if $x \in D$, then $f(x) \in f(D)$. Conversely, if $f(x) \in f(D)$ then every element of the fiber of $y = f(x)$ must lie in $D$; in  	particular, $x\in D$. Therefore, \rt{fiber-invariance} is equivalent to stating that the set $\PMV(n , \cA)\cap \R_{>0}^{E_n}$ is fiber-invariant for the partition monomial map $\ph$. It turns out that a more general statement is true.

\begin{Th}\label{T:main-fiber} Let $n \ge 4$ be even. 
Consider the action of the multiplicative group $\G=\R_{>0}^{n}$ on $\R_{>0}^{E_n}$ by 
$$\lambda\cdot x=\left(\lambda_i\lambda_jx_{\{i,j\}}\, : \, \{i,j\}\in E_n\right),\quad\text{for }\ \lambda\in\G,\, x\in\R_{>0}^{E_n}.$$
Then every $\G$-invariant subset of $\R_{>0}^{E_n}$ is fiber-invariant for the partition monomial map $\ph$ in \re{ph-map}.
\end{Th}

As explained above, the following corollary implies \rt{fiber-invariance}.

\begin{Cor} Let $\cA\subseteq\cK_2$ be a set of convex bodies closed under taking positive linear combinations. 
%Consider the partition monomial map $\ph:\R^{E_n}\to\R^{\cM_n}$ as in \re{ph-map}.
Then  $\PMV(n , \cA)\cap \R_{>0}^{E_n}$ is fiber-invariant for the partition monomial map $\ph$ in \re{ph-map}.
\end{Cor}

\begin{pf} By \rt{main-fiber} it is enough to show that $\PMV(n , \cA)\cap \R_{>0}^{E_n}$ is $\G$-invariant.
This follows from the fact that the mixed volume is linear with respect to rescaling each of the bodies
by positive scalars. Indeed, 
for every pure mixed volume configuration $\PV(K)$ in $\PMV(n , \cA)$
and any $\lambda\in\G$ we have 
$$\lambda\cdot \PV(K)=\left(\lambda_i\lambda_j\V(K_i,K_j) : {\{i,j\}}\in E_n\right)=\left(\V(\lambda_iK_i,\lambda_jK_j) : {\{i,j\}}\in E_n\right)=\PV(\lambda K),$$
where $\lambda K:=(\lambda_1K_1,\dots, \lambda_nK_n)$. As $\cA$ is closed under taking positive linear combinations, we 
have $\lambda K\in\cA$ and so $\PV(\lambda K)\in\PMV(n , \cA)$.
\end{pf}

To prove \rt{main-fiber} we give an equivalent  additive version of this result by 
``linearizing'' the partition monomial map $\ph$ using the coordinate-wise logarithm.
%the restriction of $\ph$ to the positive orthant by taking the logarithm of the variables and considering the linear map $f=\ln(\phi):\R^{E_n}\to\R^{\cM_n}$. 
More precisely, let $y_e=\ln x_e$ for $e\in E_n$. Then $f=\ln\circ\,\ph\,\circ\exp$ is the linear map
\begin{equation}\label{e:f-map}
f:\R^{E_n}\to\R^{\cM_n},\quad f(y)= \left(y_M\,:\,  M\in\cM_n\right),\quad\text{where }\ y_M=\sum_{e\in M}y_e.
\end{equation}
We have the analogous action of the additive group $\Lambda=\R^{n}$ on $\R^{E_n}$:
\begin{equation}\label{e:additive-action}
\lambda\cdot x=\left(\lambda_i+\lambda_j+y_{\{i,j\}}\, : \, \{i,j\}\in E_n\right),\quad\text{for }\ \lambda\in\Lambda,\, y\in\R^{E_n}.
\end{equation}
Clearly, a subset $D\subset\R_{>0}^{E_n}$ is fiber-invariant for $\ph$ if and only if $\ln(D)\subset\R^{E_n}$ is fiber-invariant for~$f$.
Moreover, $D$ is $\G$-invariant if and only if $\ln(D)$ is $\Lambda$-invariant. Therefore, \rt{main-fiber} is equivalent to the following theorem.

\begin{Th}\label{T:additive-fiber} Let $n \ge 4$ be even. 
Consider the action of the additive group $\Lambda=\R^{n}$ on $\R^{E_n}$ in \re{additive-action}.
Then every $\Lambda$-invariant subset of $\R^{E_n}$ is fiber-invariant for the linear map $f$ in \re{f-map}.
\end{Th}

\begin{pf} Let $g:\R^n\to \R^{E_n}$ be a linear map given by $g(\lambda)=\left(\lambda_i+\lambda_j\, : \, \{i,j\}\in E_n\right)$.
Then $D\subset\R^{E_n}$ is $\Lambda$-invariant if and only if $D+\im g=D$. Next, observe that 
$f(\1)=\frac{|E_n|}{2}\1$, where $\1$ denotes the vector of all 1's in the corresponding space. In other words, 
$\im f$ contains the line $\la\1\ra$ spanned by~$\1$. Thus, we obtain a sequence
\begin{equation}\label{e:exact-seq}
\begin{tikzcd}   \R^{n} \arrow[r, "g"] & \R^{E_n} \arrow[r, "\bar f"] & \R^{\cM_n}/\la\1\ra\,, \end{tikzcd}
\end{equation}
where $\bar f(y)=f(y)+\la\1\ra$. 

Let us check that this sequence is exact in $\R^{E_n}$.
First, the image of a standard basis vector $e_i\in\R^n$ under $g$ is the vector $v=(v_e)_{e\in E_n}$ where
$v_e=1$ if $i\in e$ and $0$, otherwise. As every matching $M\in\cM_n$ has precisely one edge connected to vertex $i$, we see that $f( g(e_i))=f(v)=\1$. Thus, $\bar f\circ g=0$. On the other hand, Corollary~2.3 in \cite{lovasz1987matching} implies that 
$\rk\, f=|E_n|-n+1$. (The matrix for $f$ in the standard bases is the incidence matrix of perfect matchings in
a complete graph.) Hence, $\rk\, \bar f=|E_n|-n$. Also, it is easy to see that $\rk\, g=n$.
This shows that the sequence \re{exact-seq} is exact in~$\R^{E_n}$.
By \rl{exact-seq}, $D$ is fiber-invariant for the linear map $\bar f$ and, hence, for $f$ as well.
\end{pf}

%%%%%%%%%%%%%%%%%%

\section{Rigid configurations and polyhedrality cases of $C_{n,d}$}\label{S:rigid}

In this section we study so called ``rigid'' configurations of points in a projective space, which allow us to compute $C_{n,d}$ in the special cases 
$(n,d)\in\{(4,2),(6,2),(6,3)\}$. 
%In this section we state our results for $(n,d)=(6,2)$ in the centrally-symmetric case and we show non-polyhedrality for $(2k, 2)$ when $k \ge 4$.  
%Recall that the class of planar centrally-symmetric bodies coincides with $\cZ_2$. 

Let $A$ be a configuration (i.e. a multiset) of points in a projective space $\RP^\ell$. By its {\it span} $\spn(A)$ we mean the smallest projective subspace of $\RP^\ell$ containing $A$. We say a configuration $A$ is {\it non-degenerate} if $\spn(A)=\RP^\ell$.
Denote by $\cH(A)$ the set of all projective hyperplanes of the form $\spn(A')$ for some $A'\subset A$.

\begin{Def} Let $A$ be a non-degenerate point configuration in $\RP^\ell$. We say a point $a\in A$ is {\it locked}
if it lies on at least $\ell$ distinct hyperplanes in $\cH(A \setminus \{a\})$, and is called {\it free} otherwise. A configuration $A$ is called \emph{rigid} if every point in 
$A$ is locked.
\end{Def}

It is easy to see that any point configuration in $\RP^\ell$ whose points have multiplicity greater than one is rigid.
Clearly, for $\ell=1$ every rigid configuration has this form. Expectedly, the case of $\ell=2$ is more interesting. The next proposition describes all rigid configurations of size up to 6 in $\RP^2$, up to projective transformations.

%In the case of $d=3$, we consider non-degenerate (i.e. non-collinear)
% point configurations $A$ in $\RP^2$. We let 
%Point configurations  all of whose points are locked are called rigid. 
%Let $\cA$ be a collection of $n$ points in $\RP^2$ (not necessarily distinct, so $\cA$ is a multiset). A point 
%$a\in\cA$ is called {\it simple} if it has multiplicity one, otherwise it is called {\it multiple}.
%We will assume that $\cA$ is non-trivial, that is the points are not all collinear. 
%Let $L_{\cA}$ be the set of distinct lines connecting every pair of distinct elements in $\cA$.
%We say $a\in\cA$ is {\it locked} if it is incident to at least two lines in 
%$L_{\cA\setminus\{a\}}$, otherwise it is called {\it free}.
%We say  $\cA$ forms a {\it rigid configuration} if every $a\in\cA$ is locked. 
%Clearly, if every point in $\cA$ is multiple then  $\cA$ is rigid. 
%We have the following proposition.

\begin{Prop}\label{P:rigid} 
There are no rigid configuration with fewer than 6 points in $\RP^2$. There are exactly two rigid configurations of 6 points in $\RP^2$, up to projective transformations. They are depicted in \rf{6-rigid}.
\end{Prop}

\begin{figure}[h]
\begin{center}
\includegraphics[scale=.39]{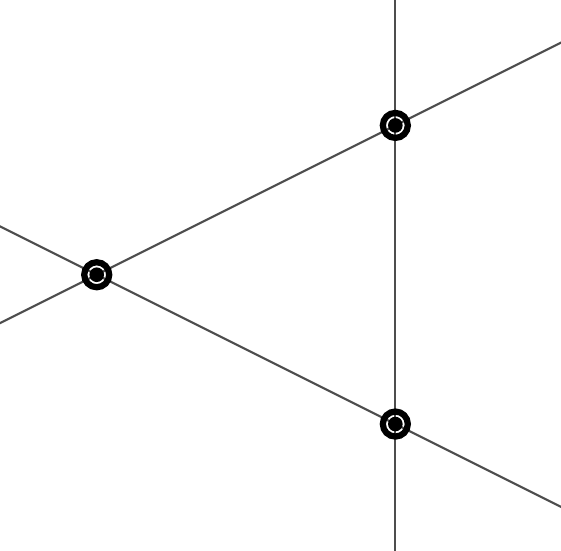}\hspace{1.5cm}
\includegraphics[scale=.43]{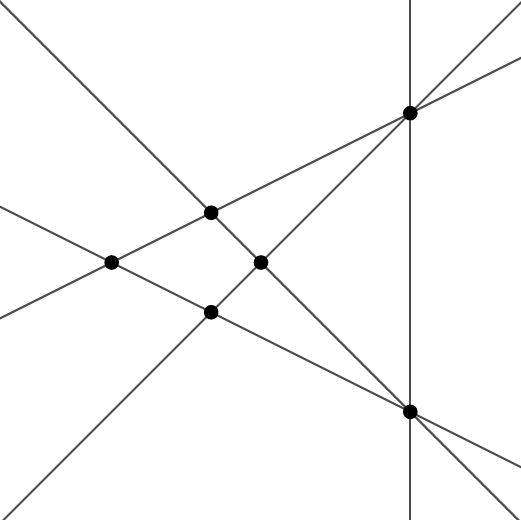}
\end{center}
\caption{Two rigid configurations with six points, counting multiplicities}
\label{F:6-rigid}
\end{figure}

\begin{pf} Let $A$ be a non-degenerate rigid configuration with up to 6 points. We let $|A|$ denote the number of elements of $A$
(counting multiplicities) and $n_A$ denote the number of %corresponding set of 
distinct elements of $A$. By definition, $\cH(A)$ is the set of distinct lines connecting every pair of distinct elements in $A$ and $a\in A$ is locked if and only if it is incident to at least two lines in $\cH(A\setminus\{a\})$. 

As $A$ is non-degenerate, $n_A\geq 3$. If $n_A=3$ then $|\cH(A)|=3$ and, hence, $A$ is rigid only
if every point in $A$ is multiple. When $|A|=6$ this produces one such configuration with
three non-collinear points of multiplicity two each.

Suppose $n_A=4$ and pick $a\in A$ such that $A\setminus\{a\}$ is non-degenerate and, hence, $|\cH({A\setminus\{a\}})|=3$. We may assume that $a$ is simple, as $A$ has at most three distinct collinear points and at most two multiple points. But 
if $a$ is locked, it must lie on two of the lines in $\cH({A\setminus\{a\}})$, i.e. coincide with a point
in $A\setminus\{a\}$, a contradiction.

Suppose $n_A=5$ and let $a,b\in A$ such that $A\setminus\{a,b\}$ is non-degenerate. As above, 
we may assume that $a$ and $b$ are simple points.
Assume $a$ lies on one of the three lines in $\cH({A\setminus\{a,b\}})$. Then $\cH({A\setminus\{b\}})$
consists of four lines (see \rf{four-points}). 
\begin{figure}[h]
\begin{center}
\includegraphics[scale=.42]{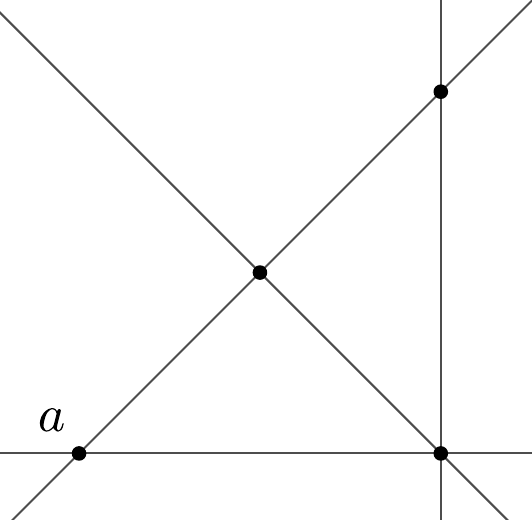}
\end{center}
\caption{The four lines in $\cH({S\setminus\{b\}})$}
\label{F:four-points}
\end{figure}
If $b$ is locked it must lie on at least two of the four lines, which is 
impossible if $n_A=5$. Assume $a$ does not lie on any of the lines in $\cH({A\setminus\{a,b\}})$. Then
$\cH({A\setminus\{b\}})$ consists of six lines with seven incidence points: the four points in $A\setminus\{b\}$ and
three extra points (see \rf{seven-points}). 
\begin{figure}[h]
\begin{center}
\includegraphics[scale=.42]{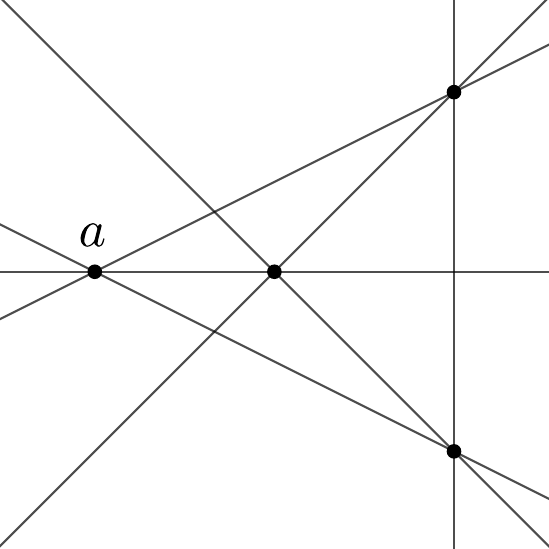}
\end{center}
\caption{The six lines in $\cH({A\setminus\{b\}})$}
\label{F:seven-points}
\end{figure}
If $b$ is locked it must be one of the extra points. But then
$a$ lies on only one line in $\cH({A\setminus\{a\}})$, i.e. is free, a contradiction.

Finally, suppose $n_A=6$, so all points in $A$ are simple. Fix $a_1\in A$. Since $a_1$ is locked, there exist 
two disjoint 2-element subsets $\{b_1, b_2\}$ and $\{c_1,c_2\}$ of $A$ such that
$\{a_1,b_1,b_2\}$ and $\{a_1,c_1,c_2\}$ are triples of collinear points. Let $a_2$ be the remaining point in $A$.
Since $b_1$ is locked, it must appear in another collinear triple, disjoint from $\{a_1,b_2\}$. We may assume this 
triple is $\{a_2,b_1,c_1\}$ (up to ordering of $c_1,c_2$). Similarly, $\{a_2,b_2,c_2\}$ must be another collinear triple. This corresponds to exactly one rigid configuration in \rf{6-rigid}, up to projective transformations. 
\end{pf}

\begin{Rem} There exist infinite families of rigid configurations with $n\geq 7$. For example, \rf{7-rigid-points}
represents a 1-parameter family of rigid configurations parametrized by the slope of the dotted line.
\end{Rem}

\begin{figure}[h]
\begin{center}
\includegraphics[scale=.42]{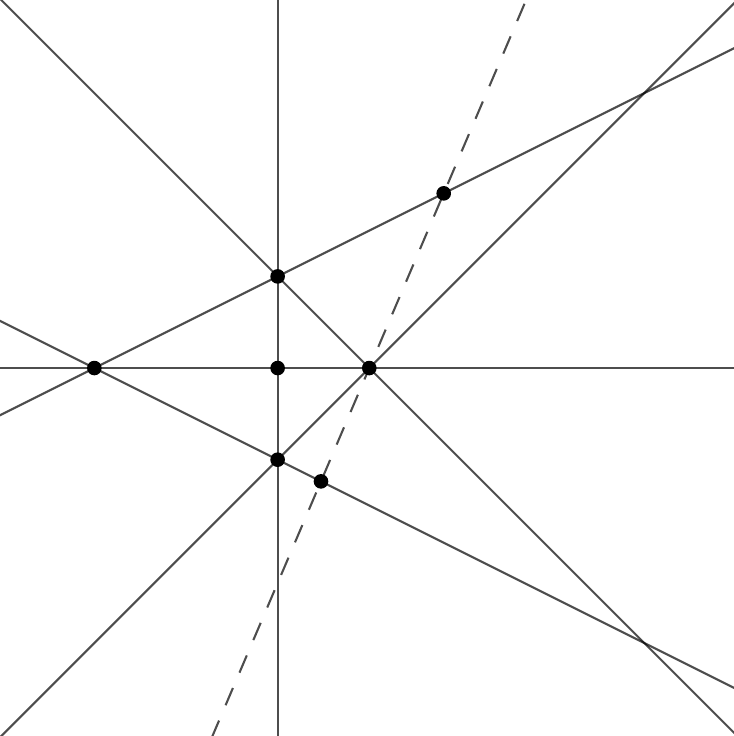}
\end{center}
\caption{A 1-parameter family of rigid configurations of size 7}
\label{F:7-rigid-points}
\end{figure}

\begin{Th}\label{T:conic-hull}  
Consider the map $\Phi$ in \re{our-map} and let $C_{n,d}$ be the conic hull of $\Phi(\cZ_d^n)$.
Assume that either $(n,d)=(6,3)$ or
$(n,d)=(2k,2)$ for some $k\geq 2$. Then $C_{n,d}$ is 
generated by the values $\Phi(L)$ over all $n$-tuples of segments $L$ whose corresponding point configurations in $\RP^{d-1}$ are rigid. Moreover, $C_{4,2}$, $C_{6,2}$, and $C_{6,3}$ are polyhedral.
\end{Th}

\begin{pf} %Throughout the proof we use the identifications and notation introduced in \rs{proj-view}. 
Let $\cL_d^n$ be the space of all $n$-tuples $L=(L_1,\dots, L_n)$, where each $L_i$ is a segment in $\R^d$.  
Recall that each zonoid is a Hausdorff limit of zonotopes (i.e. finite sums of segments) and that $\Phi$ is Minkowski multilinear and continuous with respect to the Hausdorff metric. This implies that $C_{n,d}$ is the conic hull of $\{\Phi(L) :  L\in\cL_d^n\}$. Thus, it is enough to prove that the extremal rays of $C_{n,d}$ come from rigid point configurations
in $\RP^{d-1}$.

Using our identification from \rs{proj-view}, we write $\Phi(U)$ instead of $\Phi(L)$ where $U=(u_1,\dots, u_n)$ is the $n$-tuple of vectors corresponding to $L\in\cL_d^n$ and we write $A=\{a_1,\dots, a_n\}$ for the corresponding point configuration in $\RP^{d-1}$.

%\noindent 
{\it Case $d=2$.} First we look at the easier case of $d=2$ and $n=2k$ for some $k\geq 2$. Choose
$U=(u_1,\dots, u_n)$ with $\Phi(U)\neq 0$ and consider the corresponding point configuration 
$A=\{a_1,\dots, a_n\}$ in $\RP^1$. Without loss of generality we may assume that
$u_1,\dots, u_n$ are ordered counterclockwise, which gives us a natural circular
ordering of  $a_1,\dots,a_n$.

%Up to a projective transformation, we may assume that $A$ does not contain the infinite point, so $U=\left(\begin{matrix} a_1 & \dots & a_n\\ 1 & \dots & 1\end{matrix}\right)$, where we also assume $a_1\leq\dots\leq a_n$.

Note that the number $n_A$ of distinct elements of $A$ is at least 2 and the multiplicity of each point in $A$ is at most $k$, otherwise, $\Phi(U)=0$. 
If $A$ is not rigid then it contains a point, say $a_2$, of multiplicity one and, hence, $n_A\geq 3$.
%Again, by performing a projective transformation if needed, we may assume that $i>1$.
The idea is to move $a_2$ until it coincides with either $a_{1}$ or $a_{3}$ and, hence, no longer free. More formally, let $U(x)=(u_1,x,u_3,\dots, u_n)$ where $x\in\R^2$ is a variable vector.
%$$U(x)=\left(\begin{matrix} a_1 &  \dots & x & \dots & a_n\\ 1 & \dots & 1 & \dots & 1\end{matrix}\right),$$
%where $a_i$  is replaced with $x$, and let $U(x)$ be the corresponding $n$-tuple of vectors.
Note that each entry of $\Phi(U(x))$ is a constant multiple of $|\det(u_i,x)|$ for some $i$.
Hence, $F(x):=\Phi(U(x))$ is a continuous piecewise linear vector function of $x$, which is linear on the pointed cone generated by $u_1$ and $u_3$. 
Write $u_2=v_1+v_3$, where $v_1$ and $v_3$ are proportional to $u_1$ and $u_3$, respectively, and 
let $U'$ (resp. $U''$) be the tuple obtained from $U$ by replacing $u_2$ with $v_1$ (resp. $v_3$). Then,
by the linearity of $F(x)$ on $\cone(u_1,u_3)$, we have
\begin{equation}\label{e:pos-comb}
\Phi(U)=F(u_2)=F(v_1)+F(v_3)=\Phi(U')+\Phi(U'').
\end{equation}
Note that the corresponding point configurations $A'$ and $A''$ have fewer free points than $A$.
%for some positive $\lambda,\mu$. If $F(a_{i-1})$ and $F(a_{i+1})$ span two different rays then \re{pos-comb} shows that the ray spanned by $\Phi(U)$ is not extremal. Otherwise, we
%may replace $U$ with either $U(a_{i-1})$ or $U(a_{i+1})$, which now has fewer free points and generates the same
%ray as $U$. 
%Continuing this way we show that either the ray spanned by $\Phi(U)$ is non-extremal or comes from a rigid configuration.
Continuing this way we can express every non-zero vector $\Phi(U)$  in $C_{n,2}$ as a sum of vectors $\Phi(U')$ whose corresponding point configurations $A'$ have no free points, i.e. are rigid configurations. It remains to note that for $n=4$ (resp. $n=6$) there is a unique, up to a projective transformation, configuration of two (resp. three) double points in $\RP^{1}$. This implies that $C_{4,2}$
and $C_{6,2}$ are polyhedral.

{\it Case $(n,d)=(6,3)$.} We use a similar idea of gradually discarding the configurations that are redundant in generating $C_{6,3}$ and coming up with a generating set of rigid configurations that correspond to finitely many rays of $C_{6,3}$. 

Consider a tuple of vectors $U=(u_1,\dots,u_6)$ in $\R^3$ such that $\Phi(U)\neq 0$ and let $A=\{a_1,\dots, a_6\}$ be the corresponding non-degenerate point configuration in $\RP^2$. Note that $A$ may have at most double points, as otherwise $\Phi(U)= 0$.

%As before, we may assume that $A$ does not contain points on the infinite line, so $u_i=(a_i,1)\in\R^3$ for $1\leq i\leq 6$.

{\it Step 1.} Here is the first step of reduction. Assume without loss of generality that $\{u_1,\ldots, u_5\}$ spans $\R^3$ and replace $u_6$ by a varying vector $x\in\R^3$. Denote $U(x)=(u_1,\dots, u_5,x)$. 
Then $F(x):=\Phi(U(x))$ is a continuous piecewise linear vector function, as the components of $F(x)$, up to a constant factor, have the form $|\det(u_i,u_j,x)|$. Consider the central plane arrangement $\cH$ of the planes $\{ x \in \R^3 \colon \det(u_i,u_j,x) =0 \}$ defined by all $\{i,j\} \subseteq [5]$ for which $u_i$ and $u_j$ are linearly independent. Since  $\{u_1,\ldots, u_5\}$ spans $\R^3$, the arrangement $\cH$ is pointed, i.e. each maximal cell of the arrangement is a pointed three-dimensional cone. Pick a cell $R$ of the arrangement which contains the vector $u_6$. As each determinant $\det(u_i,u_j,x)$ does not change sign
as $x$ varies over $R$, the function $F(x)$  is linear on $R$. Pick vectors $v_1,\dots, v_r$ generating the extremal rays of $R$, such that $u_6 = v_1 + \cdots + v_r$. By the linearity of $F$ on $R$, we have %$F(x) = F(v_1) + \cdots + F(v_m)$.
%This means that 
$\Phi(U) = \Phi(U_1) + \cdots + \Phi(U_r)$, where $U_s = (u_1,\dots, u_5,v_s)$ for $1\leq s\leq r$.

Let us analyze the point configuration $A_s$ corresponding to $U_s$. Denote by $b_s$ the point in
$\RP^2$, corresponding to $v_s$. Note that $b_s$ lies on the two distinct lines with homogeneous
equations $\det(u_i,u_j,x) =0$ and $\det(u_\ell,u_m,x) =0$, for some $\{i,j\},\{\ell,m\}\subset[5]$.
If $a_i,a_j,a_\ell,a_m,b_s$ are all distinct then $A_s$ contains two triples of collinear points 
$\{a_i,a_j,b_s\}$ and $\{a_\ell,a_m,b_s\}$, see the left diagram in \rf{Step1}. Otherwise, $b_s$
coincides with one of the $a_i$, $a_j$, $a_\ell$, $a_m$, as depicted in the right diagram in \rf{Step1}.
\begin{figure}[h]
\begin{center}
\tikz[scale=0.6,baseline=2mm]{\draw[thick] (0,2) -- (0,0) -- (2,0); \fill (0,0) circle (0.15) (0,1) circle (0.15) (0,2) circle (0.15) (1,0) circle (0.15) (2,0) circle (0.15);}\hspace{3cm}
\tikz[scale=0.6,baseline=2mm]{\draw[thick] (0,2) -- (0,0) -- (2,0); \fill (0,0) circle (0.15) (0,1) circle (0.15) (0,2) circle (0.15)  (2,0) circle (0.15); \draw (0,0) circle (0.25);}
\end{center}
\caption{Two types of special configurations in Step 1.}
\label{F:Step1}
\end{figure}
To summarize, to generate $C_{6,3}$ it is enough to use non-degenerate configurations containing either (a)
two adjacent collinearities or (b) a double point collinear with two more points. 

{\it Step 2(a).} Let $A$ be a non-degenerate configuration without multiple points which has adjacent collinearities, say, $\{a_1,a_2,a_3\}$ and $\{a_1,a_4,a_5\}$. As in Step 1, let $U(x)=(u_1,\dots, u_5,x)$
and $F(x) = \Phi(U(x))$ be the corresponding piecewise linear function. Consider the pointed plane arrangement $\cH$ with planes $\{ x\in\R^3 : \det(u_i,u_j,x) =0\}$, for $\{i,j\}\subset[5]$ such that  $u_i$ and $u_j$ are linearly independent. Let $R$ be the cell of $\cH$ containing $u_6$ and write
$u_6 = v_1 + \cdots + v_r$ for some generators $v_1,\dots, v_r$ of extremal rays of $R$. As above, by
the linearity of $F(x)$ on $R$, we have
$\Phi(U) = \Phi(U_1) + \cdots + \Phi(U_r)$, where $U_s = (u_1,\dots, u_5,v_s)$ for $1\leq s\leq r$.
Let $A_s=\{a_1,\dots, a_5,b_s\}\subset\RP^2$ be the corresponding point configuration.
As in Step 1, either $A_s$ has a double point collinear to a pair of other points or
$A_s$ has no multiple points and two triples of collinear points involving $b_s$. In this case
$A_s$ is a rigid configuration as in the right of \rf{6-rigid}.

{\it Step 2(b).} Now we take care of the configurations with a double point collinear with two more points. Assume $a_1 = a_5$ and $\{a_1,a_2,a_3\}$ belong to the same line. Note that neither 
$a_4$ nor $a_6$ can lie on this line, as otherwise, $\Phi(U)= 0$. Let $\cH$ be the arrangement of planes $H_i=\{ x\in\R^3  : \det(u_1,u_i,x) = 0\}$, for $i=2,3,4$. Note that
$H_2=H_3$, since $\{u_1,u_2,u_3\}$ are coplanar. Also $H_3\cap H_4$ is
the line spanned by $u_1$. Therefore, $\cH$ has four maximal cells, each being the 
direct product of a line and a two-dimensional pointed cone. As before, the function
$F(x)=\Phi(U(x))$, where $U(x)=(u_1,\dots, u_5, x)$, is continuous and linear on the cells of $\cH$.
Let $R$ be the cell containing $u_6$ and write $u_6=v_3+v_4$ for $v_3,v_4$ in the facets of $R$,
such that $v_i\in H_i$, for $i=3,4$.
As before, since $F$ is linear on $R$, $\Phi(U)=\Phi(U_3)+\Phi(U_4)$ where $U_i=(u_1,\dots, u_5,v_i)$, for $i=3,4$. Let $A_i=(a_1,\dots, a_5,b_i)$ be the corresponding point configurations, for $i=3,4$. Observe that the configuration $A_3$ 
now has collinearity $\{a_1,a_2,a_3,b_3\}$ and, hence, $\Phi(U_3)=0$. The configuration
$A_4$ consists of two adjacent collinearities $\{a_1,a_2,a_3\}$ and $\{a_1,a_4,b_4\}$ where
$a_1=a_5$ is a double point. We have thus shown that each configuration containing  a double point collinear with two more points can be replaced with the configuration as in \rf{Step2}, without changing
the span of $\Phi(U)$.
\begin{figure}[h]
\begin{center}
\tikz[scale=0.6,baseline=2mm]{\draw[thick] (0,2) -- (0,0) -- (2,0); \fill (0,0) circle (0.15) (0,1) circle (0.15) (0,2) circle (0.15) (1,0) circle (0.15) (2,0) circle (0.15); \draw (0,0) circle (0.25);}
\end{center}
\caption{Special configuration in Step 2(b).}
\label{F:Step2}
\end{figure}

{\it Step 3.} The final step is to show that the configuration as in \rf{Step2} can be further reduced
to the rigid configuration consisting of three double points (see \rf{6-rigid}) without changing the
span of $\Phi(U)$. So let $a_1=a_6$ and $\{a_1,a_2,a_3\}$, $\{a_1,a_4,a_5\}$ be two
collinear triples. Note that adding multiples of $u_1$ to $u_3$ and to $u_5$ does not change the span
of $\Phi(U)$. We can always add appropriate multiples of $u_1$ to $u_3$ and $u_5$ such that
$u_3$ becomes proportional to $u_2$ and $u_5$ becomes proportional to $u_4$. Projectively, 
this means that we can make $a_2=a_3$ and $a_4=a_5$ without changing the ray generated by
$\Phi(U)$, which concludes the proof.

Here is the diagram of the reduction that has been carried out above. 

\begin{center} 
\begin{tikzpicture}[
    box/.style={
        rectangle, 
        rounded corners, 
        draw=black, 
        thick, 
        minimum height=1cm, 
        minimum width=2cm, 
        align=center, 
        text centered
    },
    arrow/.style={
        ->, 
        thick, 
        draw=black
    }
]

% First box
\node[box] at (0,0) (anything) {$A$ is non-degenerate};

% Second box, positioned to the right of the first box and lower
\node[box] at (5,-1) (corner) {$A$ contains  \tikz[scale=0.3,baseline=2mm]{\draw[thick] (0,2) -- (0,0) -- (2,0); \fill (0,0) circle (0.2) (0,1) circle (0.2) (0,2) circle (0.2) (1,0) circle (0.2) (2,0) circle (0.2);}};

\node[box,fill=green!20!white] at (5,-3) (rigidcorner) {$A$ is  \tikz[scale=0.3,baseline=2mm]{\draw[thick] (1,0) -- (0,2) -- (0,0) -- (2,0) -- (0,1); \fill (0,0) circle (0.2) (0,1) circle (0.2) (0,2) circle (0.2) (1,0) circle (0.2) (2,0) circle (0.2) (2/3,2/3) circle (0.2);}};

%\node[box] at (0,-2) (doublepoint) {$A$ contains a double point \tikz[scale=0.3,baseline=-1mm]{\fill (0,0) circle (0.2) ; \draw (0,0) circle (0.4);}};

\node[box] at (0,-3) (double1) {$A$ contains \tikz[scale=0.3,baseline=-1mm]{\fill (0,0) circle (0.2)  (1,0) circle (0.2) (2,0) circle (0.2); \draw (0,0) circle (0.4) ; \draw (0,0) -- (2,0);}};

\node[box] at (0,-5) (double2) {$A$ is \tikz[scale=0.3,baseline=2mm]{\fill (0,0) circle (0.2)  (1,0) circle (0.2) (2,0) circle (0.2) (0,1) circle (0.2) (0,2) circle (0.2); \draw (0,0) circle (0.4) ; \draw (0,2) -- (0,0) -- (2,0);}};

\node[box,fill=green!20!white] at (0,-7) (threedoublepoints) {$A$ is  \tikz[scale=0.3,baseline=1mm]{\fill (0,0) circle (0.2) (1,0) circle (0.2) (0,1) circle (0.2); \draw (0,0) circle (0.4) (1,0) circle (0.4) (0,1) circle (0.4);  \draw (0,0) -- (1,0) -- (0,1) -- cycle;}};

% Draw an arrow from the first box to the second box
\draw[arrow]  (anything) -- (corner);
\draw[arrow] (corner) -- (rigidcorner);
%\draw[arrow] (anything) -- (doublepoint);
\draw[arrow] (corner) -- (double1);
\draw[arrow] (anything) -- (double1);;
\draw[arrow] (double1) -- (double2);
\draw[arrow] (double2) -- (threedoublepoints);
\end{tikzpicture} 
\end{center} 

\end{pf}

\section{Four bodies in $\R^2$: description of $\PMV(4,\cK_2)$}\label{S:framework-4-2} 
In this section we illustrate our approach in
the smallest non-trivial example $(n,d)=(4,2)$. 
Each $4$-tuple of convex bodies $K=(K_1,\dots,K_4)$ defines 
six values of mixed volumes $\V_{i j}(K)=\V(K_i,K_j)$ for $\{i,j\}\subset[4]$. 
For simplicity, we will abbreviate  $\V_{i j}=\V_{i j}(K)$ in the rest of the section.
There are three partitions of $[4]$ into 2-element subsets 
%$1\,2\,|\,3\,4$, $1\,3\,|\,2\,4$, and $1\,4\,|\,2\,3$, 
$\Prt(4,2)=\left\{12\,|\,34, \ 13\,|\,24,\ 14\,|\,23\right\}$, 
so the map $\Phi$ in \re{our-map} becomes
$$\Phi:\cK_2^4\to \R^3, \quad \Phi(K)=(\V_{12}\!\V_{34}\,, \V_{13}\!\V_{24}\,, \V_{14}\!\V_{23}).$$
Now for a $4$-tuple of line segments $L=(L_1,\dots,L_4)$ with $L_i=[0,u_i]$, let $x_{i j}=\det(u_i,u_j)$ be the corresponding Plücker coordinates. Then $\V_{ij}=|x_{ij}|$. The Plücker embedding of $\Gr_{2,4}$ into $\RP^5$
is given by the Plücker relation
$$x_{1 2}x_{3 4}-x_{1 3}x_{2 4}+x_{1 4}x_{2 3}=0.$$
This, with the help of the triangle inequality for the absolute value, produces three
 inequalities for the points of $\Abs(\Gr_{2,4})$ of the form \re{framework}:
\begin{align}\label{e:plucker}
\V_{12}\!\V_{34}\,-\,\V_{13}\!\V_{24}\,+\,\V_{14}\!\V_{23}&\geq 0\nonumber\\
\V_{12}\!\V_{34}\,+\,\V_{13}\!\V_{24}\,-\,\V_{14}\!\V_{23}&\geq 0\\
-\V_{12}\!\V_{34}\,+\,\V_{13}\!\V_{24}\,+\,\V_{14}\!\V_{23}&\geq 0.\nonumber
\end{align}
As the above inequalities are linear in each $L_i$, they also hold for the points of $\PMV(4,\cZ_2)$. 
In fact, as shown in \cite{AS23}, the inequalities \re{plucker} also hold (and
completely describe) the space $\PMV(4,\cK_2)$. 

Observe that the inequalities \re{plucker} define a polyhedral cone in $\R^3$ with coordinates
$y_{ij\,|\,kl}=\V_{ij}\!\V_{kl}=|x_{ij}x_{kl}|$ for the three partitions $ij\,|\,kl=[4]$. 
This is the cone $C_{4,2}$ as in \rt{conic-hull}. It
has three generators $(0,1,1)$, $(1,0,1)$, and $(1,1,0)$ which are the images under $\Phi$
of three quadruples of vectors 
$U_1=\left(\begin{matrix}0 & 0 & 1 & 1\\1 & 1 & 0 & 0\end{matrix}\right)$,
$U_2=\left(\begin{matrix}0 & 1 & 0 & 1\\1 & 0 & 1 & 0\end{matrix}\right)$, and
$U_3=\left(\begin{matrix}0 & 1 & 1 & 0\\1 & 0 & 0 & 1\end{matrix}\right)$.
All three quadruples correspond to the unique (up to projective transformations) 
rigid configuration $2\cdot\{0,\infty\}$ in $\RP^1$, where of $0=(0:1)$ and $\infty=(1:0)$ in projective coordinates.

\section{Six centrally-symmetric bodies in $\R^2$: description of $C_{6,2}$}\label{S:framework-6-2} 

%In this section we state our results for $(n,d)=(6,2)$ in the centrally-symmetric case and we show non-polyhedrality for $(2k, 2)$ when $k \ge 4$.  
%Recall that the class of planar centrally-symmetric bodies coincides with $\cZ_2$. 
%In fact, in dimension two, the class of zonoids $\cZ_2$ coincides with the class of all centrally-symmetric bodies.
For $(n,d)=(6,2)$ the map $\Phi$ from \re{our-map} becomes
\begin{equation}\label{e:our-map-(6,2)}
\Phi:\cK_2^6\to\R^{15},\quad \Phi(K)= \left(\V_{I_1}\!\V_{I_2}\!\V_{I_3}\,:\,  I_1|I_2|I_3=[6]\right).
\end{equation}
As before, we use the abbreviation  $\V_{i j}=\V_{i j}(K)$ throughout this section.

According to \rt{conic-hull}, to compute the conic hull $C_{6,2}$ of $\Phi(\cZ_2^6)$,
it is enough to evaluate $\Phi$ at $6$-tuples of segments corresponding to rigid configurations of $6$ points in $\RP^1$. There are exactly two such configurations, up to projective transformations, 
$A_{(3,3)}=3\cdot \{0,\infty\}$ and $A_{(2,2,2)}=2\cdot\{0,1, \infty\}$. 
%, where of $0=(0:1)$, $1=(1:1)$, and $\infty=(1:0)$ in projective coordinates. 
Let $L$ be a 6-tuples of segments and $U$ the corresponding 6-tuples of vectors.
In what follows we say that $L$ is of type $(3,3)$ (resp. of type (2,2,2)) if $L$ corresponds to the point configuration $A_{(3,3)}$ (resp. $A_{(2,2,2)}$). In the same way we
say that the tuple $U$ or the ray $\spn(\Phi(L))$ is of type $(3,3)$ (resp. of type (2,2,2)).

It is convenient to encode each tuple  $U$ of type $(3,3)$ by assigning three
labels from $\{1,\dots,6\}$ to $0\in A_{(3,3)}$ and the remaining labels to $\infty\in A_{(3,3)}$.
These labels represent the position of the corresponding vector in $U$. Thus,
a labeling is an ordered partition $O_U=(J_1,J_2)$ where $J_1\,|\,J_2\in \Prt(6,3)$. 
For example, $U=\left(\begin{matrix}0 & 0 & 0 & 1 & 1 & 1\\1 & 1 & 1 & 0 & 0 & 0\end{matrix}\right)$
corresponds to the labeling $O_U=(\{1, 2, 3\}, \{4,5,6\})$.
Similarly, tuples of
type $(2,2,2)$ are encoded by labelings $O_U=(J_1,J_2,J_3)$ for $J_1\,|\,J_2\,|\,J_3\in \Prt(6,2)$.
 For example, the tuple  $U=\left(\begin{matrix}0 & 0 & 1 & 1 & 1 & 1\\1 & 1 & 1 & 1 & 0 & 0\end{matrix}\right)$  corresponds to the labeling $O_U=(\{1, 2\}, \{3,4\},\{5,6\})$.

We use \emph{SageMath} \cite{sagemath} to produce the following combinatorial information about $C_{6,2}$. Note that
the symmetric group $\mathbf{S}_6$ acts on $\cZ_2^6$ by permuting the zonoids and, hence, acts on
$C_{6,2}$ as well.

\begin{Th} The conic hull $C_{6,2}$ of $\Phi(\cZ_2^6)$ is a polyhedral cone generated by 10 rays of type $(3,3)$ and 15 rays of type $(2,2,2)$. It has 975 facets which form 8 orbits under the action of $\mathbf{S}_6$. Representatives of each orbit are collected in Table \ref{table:n6d2:types}.
\end{Th}

\begin{table}[h]
\begin{center} 
\begin{tabular}{c||r|r|r|r|r|r|r|r}
	 & \multicolumn{8}{c}{Coefficients $\lambda_{I_1|I_2| I_3}$ of the eight inequality types} \\
	 \hline
	  Partition & Type~1 & Type~2 & Type~3 & Type~4 & Type~5 & Type~6 & Type~7 & Type~8 \\ \hline 
	 12$|$34$|$56 & -1 & -3 & -1 & -1 & -1 & -1 & -1 & -1 \\
	 16$|$23$|$45 & -1 & -1 & -1 & -1 & 0 & 0 & -1 & -1 \\
	 15$|$23$|$46 & 0 & -1 & -1 & -1 & 0 & 0 & 0 & -1 \\
	 14$|$23$|$56 & 0 & 1 & 1 & -1 & 0 & 0 & 1 & 1 \\
	 13$|$24$|$56 & 0 & -1 & -1 & 1 & -1 & 0 & 0 & -1 \\
	 13$|$25$|$46 & -1 & -1 & -1 & -1 & 0 & 0 & -1 & 0 \\
	 14$|$25$|$36 & 0 & -1 & -1 & 1 & 0 & 0 & 0 & -1 \\
	 15$|$24$|$36 & 1 & 3 & 3 & 1 & 1 & 0 & 2 & 2 \\
	 16$|$24$|$35 & 0 & -1 & 1 & -1 & 0 & 0 & 0 & 0 \\
	 16$|$25$|$34 & 1 & 3 & 1 & 1 & 1 & 0 & 1 & 2 \\
	 15$|$26$|$34 & 0 & 1 & 1 & 1 & 0 & 0 & 0 & 0 \\
	 14$|$26$|$35 & 0 & -1 & -1 & 1 & 0 & 0 & -1 & -1 \\
	 13$|$26$|$45 & 1 & 3 & 1 & 1 & 1 & 0 & 1 & 2 \\
	 12$|$36$|$45 & 1 & 1 & 1 & 1 & 0 & 1 & 0 & 0 \\
	 12$|$35$|$46 & 1 & 3 & 1 & 1 & 1 & 1 & 1 & 2 \\
\end{tabular}
\end{center}  

\caption{Coefficients of the inequalities  $\sum\limits_{I_1| I_2 |I_3} \lambda_{I_1| I_2 |I_3} \V_{I_1}\!\V_{I_2}\!\V_{I_3} \ge 0$,  representing the eight orbits of facets of $C_{6,2}$ under the action of $\mathbf{S}_6$.}
\label{table:n6d2:types}
\end{table} 

The above theorem produces 8 types of inequalities valid on $\Abs(\Gr_{2,6})$ and, hence, on $\PMV(6,\cZ_2)$.
Note that the inequalities of type 6 in Table \ref{table:n6d2:types} correspond to the Plücker-type inequalities
as in \re{plucker} for four of the six zonoids.

\begin{Rem}\label{R:discrepancy} 
It turns out that the inequalities of Type 3 are not satisfied in the general (non-symmetric) case.
The smallest example for which such an inequality fails is depicted in \rf{counterexample}.
\begin{figure}[h]
\begin{center}
\begin{tikzpicture}
\begin{scope}[scale=0.9] % Scale the figure down to 90%
    % Define the offset for positioning polygons in a line
    \coordinate (offset1) at (0, 0);
    \coordinate (offset2) at (2, 0);
    \coordinate (offset3) at (4, 0);
    \coordinate (offset4) at (5, 0); % Shifted one unit to the left
    \coordinate (offset5) at (8, 0);
    \coordinate (offset6) at (10, 0);

    % Draw line segments and polygons in a line with black circles at vertices and thicker lines
    
    % First segment
    \filldraw[thick] (0,0) -- (1,0);
    \filldraw[black] (0,0) circle (1.5pt);
    \filldraw[black] (1,0) circle (1.5pt);

    % Second segment
    \draw[thick, shift={(offset2)}] (0,0) -- (0,1);
    \filldraw[black, shift={(offset2)}] (0,0) circle (1.5pt);
    \filldraw[black, shift={(offset2)}] (0,1) circle (1.5pt);

    % Third segment
    \draw[thick, shift={(offset3)}] (-1,1) -- (0,0);
    \filldraw[black, shift={(offset3)}] (-1,1) circle (1.5pt);
    \filldraw[black, shift={(offset3)}] (0,0) circle (1.5pt);

    % Fourth segment (shifted one unit to the left)
    \draw[thick, shift={(offset4)}] (0,0) -- (1,1);
    \filldraw[black, shift={(offset4)}] (0,0) circle (1.5pt);
    \filldraw[black, shift={(offset4)}] (1,1) circle (1.5pt);

    % Fifth polygon (triangle)
    \draw[thick, fill=black!20!white, shift={(offset5)}] (-1,1) -- (0,0) -- (1,1) -- cycle;
    \filldraw[black, shift={(offset5)}] (-1,1) circle (1.5pt);
    \filldraw[black, shift={(offset5)}] (0,0) circle (1.5pt);
    \filldraw[black, shift={(offset5)}] (1,1) circle (1.5pt);

    % Sixth polygon (triangle)
    \draw[thick, fill=black!20!white, shift={(offset6)}] (0,0) -- (1,0) -- (0,1) -- cycle;
    \filldraw[black, shift={(offset6)}] (0,0) circle (1.5pt);
    \filldraw[black, shift={(offset6)}] (1,0) circle (1.5pt);
    \filldraw[black, shift={(offset6)}] (0,1) circle (1.5pt);
\end{scope}
\end{tikzpicture}

\end{center}
\caption{A collection of six bodies for which an inequality of Type 3 fails.}
\label{F:counterexample}
\end{figure}
Currently we do not know whether the inequalities of the other types hold in general, but our computer experiments suggest that this may be the case.
\end{Rem}

Next we turn to the question of how close $C_{6,2}$ approximates $\Phi(\cZ_2^6)$.
At the moment we can answer this only at the level of 2-dimensional faces of $C_{6,2}$.
As our computations confirm, $C_{6,2}$ is a cone over a 2-neighborly polytope. This means that every pair of rays in $C_{6,2}$ spans a 2-dimensional face. In the proposition below we show that none of the 2-faces spanned by a pair of rays of type $(2,2,2)$ contains a point of $\Phi(\cZ_2^6)$ in its relative interior. In particular, this shows that the set $\Phi(\cZ_2^6)$  is not convex. We start with a simple lemma where, as before, we identify a tuple of segments $L$ with a tuple of vectors $U$, see \rs{proj-view}. We also use the abbreviated notation \re{abbreviate}.
 
 \begin{Lem}\label{L:split} Let $U$ and $U'$ be two 6-tuples of vectors
 corresponding to configurations $A=3\cdot\{a,b\}$ and $A'=\{3\cdot a,\, 2\cdot b,\, c\}$
 for some distinct $a,b,c\in\RP^1$. Then $\Phi(U)$ and $\Phi(U')$ generate the same ray of $C_{6,2}$.
% Consider two tuples of vectors in $\R^2$
%$$
%U=\left(\begin{matrix}0 & 0 & 0 & 1 & 1 & 1\\ 1 & 1 & 1 & 1 & 1 & 1\end{matrix}\right)\quad\text{and}\quad
%U'=\left(\begin{matrix}0 & 0 & a & 1 & 1 & 1\\ 1 & 1 & 1 & 1 & 1 & 1\end{matrix}\right),
%$$
% for some $1\neq a\in\R$. Then $\Phi(U)$ and $\Phi(U')$ generate the same ray of $C_{6,2}$.
\end{Lem}

\begin{pf} Up to a projective transformation we may assume that
$a=0$, $b=1$, $c=\infty$, and 
$$
U=\left(\begin{matrix}0 & 0 & 0 & 1 & 1 & 1\\ 1 & 1 & 1 & 1 & 1 & 1\end{matrix}\right)\quad\text{and}\quad
U'=\left(\begin{matrix}0 & 0 & 0 & 1 & 1 & 1\\ 1 & 1 & 1 & 1 & 1 & 0\end{matrix}\right).
$$
%Assume $a\neq\infty$. Writing the labeled points of $A$ and $A'$ as column vectors we obtain
Note that the only partitions $P\in \Prt(6,2)$ for which $\V_P(U)$ and $\V_P(U')$ are non-zero are
$P=\,1\, i\, |\, 2\,j\, |\, 3\,k\,$ where $\{i,j,k\}=\{3,4,5\}$. For these $P$ we have  $\V_P(U)=\V_P(U')=1$. This shows that $\Phi(U')=\Phi(U)$ and the lemma follows. 
%The case $a=\infty$ is similar.
\end{pf}
 
 %$P=\boxed[0pt]{\,0\, i\, \big|\, 1\,j\, \big|\, 2\,k\,}$ playing with writing partitions :)

 \begin{Prop}\label{P:2-cells-(2,2,2)} Let 
 $\tau\subset C_{6,2}$ be a 2-face spanned by two rays $\rho_1,\rho_2$ of type $(2,2,2)$. Then 
 $\tau\cap \Phi(\cZ_2^6)=\rho_1\cup\rho_2$.
 \end{Prop}

 \begin{pf} Let $\Phi(U_1)$ and $\Phi(U_2)$ be generators for $\rho_1$ and $\rho_2$, respectively, and let
 $O_1,O_2$ be the corresponding labelings. Then either $O_1$ and $O_2$ share no parts or exactly one part.
 
 Assume  $O_1$ and $O_2$ share one part. By applying a projective transformation if needed, we
may assume that $O_1= (\{1,2\},\{3,4\}, \{ 5,6\})$ and $O_2=(\{1,3\},\{2,4\}, \{ 5,6\})$.
Suppose $\Phi(\cZ_2^6)$ contains a point $x=s\,\Phi(U_1)+t\,\Phi(U_2)\in\tau$, for some positive $s,t\in\R$. Then there exists $Z=(Z_1,\dots, Z_6)\in\cZ_2^6$ such that $\Phi(Z)=x$. In what follows we use $\V_I$ to abbreviate the mixed volume $\V_I(Z)$ for $I\in {[6]\choose 2}$.

First, note that  $\V_{5 \,6}=0$ and $\V_{i\,5}$, $\V_{i\,6}$ are positive for $1\leq i\leq 4$.
Hence, after a linear transformation, we have $Z_5=Z_6=[0,e_1]$. We can also rescale $Z_1,\dots, Z_4$ such that $\V_{i\,5}=\V_{i\,6}=1$ for $1\leq i\leq 4$. Computing the values of $\Phi(U_1)$ and $\Phi(U_2)$ we obtain
$$t= \V_{1\,2}= \V_{3\,4},\quad s= V_{1\,3} = V_{2\,4},\quad t+s=V_{1\,4} = V_{2\,3}.$$
On the other hand, the Plücker-type inequality for four planar bodies implies $t^2+s^2\geq (t+s)^2$, which
is impossible for positive $s,t$. Therefore, $x\in\rho_1\cup\rho_2$.

 Next, assume $O_1$ and $O_2$ do not share a part. Up to a projective transformation, we have
 $O_1= (\{1,2\},\{3,4\}, \{ 5,6\})$ and $O_2=(\{1,i\},\{3,j\}, \{ 5,k\})$,
 where $i\ne 2$, $j\ne 4$, and $k\ne 6$.  Thus, either (a) 
 $O_2=(\{1,6\},\{2,3\}, \{ 4,5\})$ or (b)
 $O_2=(\{1,4\},\{3,6\}, \{ 2,5\})$.
Again, computing the coordinates of $x=s\,\Phi(U_1)+t\,\Phi(U_2)\in\tau$, in case (a) we obtain
$$\V_{\,1\,2\,}V_{\,3\, 6\,}\V_{\,4\, 5\,}=0,\quad \V_{\,1\,2\,}\V_{\,3\, 4\,}\V_{\,5\, 6\,}=t,\quad \V_{\,1\,6\,}\V_{\,2\, 3\,}\V_{\,4\, 5\,} =s,\quad \V_{\,1\,4\,}\V_{\,2\, 5\,}\V_{\,3\, 6\,}=s+t,$$
which implies that either $s=0$ or $t=0$ (or both). Case (b) is similar.
 \end{pf}

It turns out that for the 2-faces spanned by two rays of type (3,3) the situation is opposite: such faces
are part of  $\Phi(\cZ_2^6)$. 

 \begin{Prop}\label{P:2-cells-(3,3)} Let 
 $\tau\subset C_{6,2}$ be a 2-face spanned by two rays of type $(3,3)$. Then 
 $\tau\subset \Phi(\cZ_2^6)$.
 \end{Prop}

\begin{pf} %Applying a projective transformation, we may replace $A_{(3,3)}$ with $A=3\cdot\{0,1\}$.
% consists of $(0:1)$ and $(1:1)$, each with multiplicity 3. 
Let $\Phi(U_1)$ and $\Phi(U_2)$ be generators for $\rho_1$ and $\rho_2$, respectively, and let $O_1,O_2$ be the corresponding labelings. Up to a projective transformation, we may assume that 
 $O_1= (\{1,2, 3\},\{4, 5,6\})$ and $O_2=(J,J^c)$, where $J=\{1,2,i\}$ for some $i\in\{4,5,6\}$. 
We assume $i=4$, so  $O_2= (\{1,2,4\},\{3,5,6\})$. The other cases are completely analogous. 
By \rl{split}, we can replace $A_{(3,3)}$ with $A_1=\{2\cdot 0,\, 1,\, 3\cdot\infty\}$ and
$O_1$ with  $O_1'= (\{1,2\}, \{3\},\{4, 5,6\})$ without changing the ray $\rho_1$.
Similarly, we can replace $A_{(3,3)}$ with $A_2=\{3\cdot 0,\, 1,\, 2\cdot\infty\}$ and
$O_2$ with  $O_2'= (\{1,2,4\}, \{3\},\{5,6\})$ without changing the ray $\rho_2$.
This means that we may replace $U_1$ and $U_2$ with
$$
U_1'=\left(\begin{matrix}0 & 0 & 1 & 1 & 1 & 1\\ 1 & 1 & 1 & 0 & 0 & 0\end{matrix}\right)
\text{ and }
U_2'=\left(\begin{matrix}0 & 0 & 1 & 0 & 1 & 1\\ 1 & 1 & 1 & 1 & 0 & 0\end{matrix}\right),
$$
respectively, without changing the rays $\rho_1$ and $\rho_2$.

Now consider a point $x=s\,\Phi(U_1')+t\,\Phi(U_2')\in\tau$, for some $s,t\in\R_{\geq 0}$. Define $Z_1=Z_2=[0,e_2]$,
$Z_3=[0,e_1+e_2]$, $Z_4=s[0,e_1]+t[0,e_2]$, and $Z_5=Z_6=[0,e_1]$. (Here $Z_4$ is a parallelogram and the rest are segments.) Then, by linearity, $\V_P(Z)=s\V_P(U_1')+t\V_P(U_2')$, for every $P\in \Prt(6,2)$. This shows that $x=\Phi(Z)$, i.e., $x\in \Phi(\cZ_2^6)$.
\end{pf}

For 2-faces generated by rays of different type the situation is more subtle. Assume
$\rho_1$  has type $(3,3)$ and $\rho_2$ has type $(2,2,2)$ and let $O_1=(J_1,J_2)$ and $O_2=(I_1,I_2,I_3)$ be the corresponding labelings.
We write $O_2\preccurlyeq O_1$ if $I_i\subset J_j$ for some $i\in[3]$ and $j\in[2]$.

 \begin{Prop}\label{P:(2,2,2)-(3,3)} Let 
 $\tau\subset C_{6,2}$ be a 2-face spanned by two rays $\rho_1,\rho_2$ of types $(3,3)$
 and $(2,2,2)$ with labelings $O_1$ and $O_2$, respectively.
If $O_2\preccurlyeq O_1$ then $\tau\subset \Phi(\cZ_2^6)$. Otherwise, $\tau\cap \Phi(\cZ_2^6)=\rho_1\cup\rho_2$.
 \end{Prop}

 \begin{pf} 
 %Let $A_1$ and $A_2$ be the configurations corresponding to $\rho_1$ and $\rho_2$ multilabeled by $P_1\in\Prt(2,2,2)$ and $P_2\in\Prt(3,3)$, respectively. Without loss of generality we may assume that $A_1=(0,1,\infty)$ and $A_2=(0,\infty)$.
 Let $\Phi(U_1)$ and $\Phi(U_2)$ be generators for $\rho_1$ and $\rho_2$, respectively, and let $O_1,O_2$ be the corresponding labelings.
 
Suppose $O_2\preccurlyeq O_1$. By applying a projective transformation if needed, we
 may assume that $O_1=(\{1,2,3\},\{4,5,6\})$ and $O_2=(\{1,2\},\{3,4\},\{5,6\})$. By \rl{split} we can replace $A_{(3,3)}$ with $A=\{3\cdot 0,1,2\cdot\infty\}$ labeled by $O_2'=(\{1,2,3\},\{4\},\{5,6\})$ without changing the ray $\rho_2$.
Now to show that $\tau\subset \Phi(\cZ_2^6)$ it is enough to show that for any $\lambda\in[0,1]$ we have $(1-\lambda)\Phi(U_1)+\lambda \Phi(U_2')=\Phi(Z)$ for some $Z\in\cZ_2^6$. In fact, we can take $Z$ to be a collection of segments:
 $Z_1=Z_2=[0,e_2]$, $Z_3=[0,(1-\lambda)e_1+\lambda e_2]$, $Z_4=[0,e_1+e_2]$, and $Z_5=Z_6=[0,e_1]$. 
  
 If $O_2\not\preccurlyeq O_1$ then, up to permuting the labels and the points, we have
 $O_1=(\{1,3,5\},\{2,4,6\})$ and $O_2=(\{1,2\},\{3,4\},\{5,6\})$.
 Computing the coordinates of $x=s\Phi(U_1)+t\Phi(U_2)\in\tau$, for $s,t\in\R_{\geq 0}$, we obtain
$$\V_{1\,3}\V_{2\, 4}\V_{5\, 6}=0,\quad \V_{1\,2}\V_{3\, 4}\V_{5\, 6}=t,\quad \V_{1\,3}\V_{2\, 5}\V_{4\, 6}=s,\quad \V_{1\,5}\V_{2\,4}\V_{3\, 6}=s,$$
which implies that either $s=0$ or $t=0$.
 \end{pf}

\section{Six zonoids in $\R^3$: description of $C_{6,3}$}\label{S:framework-6-3} 

 Let $K=(K_1,\dots, K_6)$ be a 6-tuple of convex bodies in $\R^3$. The partitions in $\Prt(6,3)$ are simply $I| I^c=[6]$
 for every 3-element subset $I\subset[6]$, where $I^c=[6]\setminus I$ denotes the complement. In this case
the map $\Phi$ in \re{our-map} is
 
 \begin{equation}\label{e:framework-3-6}
\Phi:\cK_3^6\to\R^{10},\quad \Phi(K)= \left(\V_I(K)\V_{I^c}(K)\,:\,  I\subset[6], |I|=3\right).
\end{equation}
In what follows we abbreviate $\vv_{I}=\V_I(K)\V_{I^c}(K)$. Clearly, $\vv_{I}=\vv_{I^c}$.
%Given $I,J\subset[6]$ with $|I|=3$ and $|J|=2$, we write $I>J$ if and only if $I\supset J$ or $I^c\supset J$.
 
 %\subsection{Rigid configuration of points in $\RP^{2}$}
 Recall from \rt{conic-hull} that $C_{6,3}$ %(the conic hull of $\Phi(\cZ_{2}^6)$) 
 is generated by $6$-tuples of segments corresponding  to two types of rigid configurations in $\RP^{2}$. All configurations
 of the same type are projectively equivalent, so to compute $C_{6,3}$ it is enough to consider, for example,
 $$A_1=2\cdot\{(0:0:1),(1:0:1),(0:1:1)\}$$
  (a triple of double points) and
$$A_2=\{(0:0:1),(1:0:1),(0:1:1),(1:1:1), (1:0:0),(0:1:0)\}.$$
 (six simple points), see \rf{6-rigid} for the pictures.
As confirmed by {\it SageMath}, all possible ways to label the points in $A_1$ using numbers $\{1,\dots, 6\}$ produce 15 generators of  $C_{6,3}$ of the first type. Similarly, there are $15$ more generators of  $C_{6,3}$ of the second type coming from all labelings of the points in $A_2$. We use {\it SageMath} to compute the facets of $C_{6,3}$. The facet inequalities produce 130 inequalities for the mixed volumes
of 6 zonotopes in $\R^3$. Unlike the case of $(n,d)=(6,2)$ in \rs{framework-6-2}, we can state these inequalities in an invariant way and derive them from the Grassmann-Plücker relations for the Grassmannian $\Gr_{3,6}$ and the triangle inequality, similar to the case of $(n,d)=(4,2)$ in \rs{framework-4-2}. We do this in the theorem below.

\begin{Th}  For any 6-tuple $Z=(Z_1,\dots, Z_6)$ of zonoids in $\R^3$ the quadratic monomials $\vv_{I}=\V_I(Z)\V_{I^c}(Z)$ satisfy 130 linear inequalities. Ten of them are
$\vv_I\geq 0$ and the remaining ones are 
\begin{equation}\label{e:6-zonoids-in-3d}
%0\leq \sum_{k\not\in I\supset J}\!\!\vv_{I}-\vv_{J\cup \{k}\leq \sum_{I\not > J}\vv_I,
\vv_{J\cup \{\ell\}}\leq \sum_{J\subset I\not\ni \ell}\!\vv_I\leq \vv_{J\cup \{\ell\}}+\sum_{J\not\subset I\ni \ell} \vv_I,
\end{equation}
for every choice of $J\subset[6]$ with $|J|=2$ and $\ell\in J^c$. The inequalities describe a full-dimensional pointed cone in $\R^{10}$
with 30 generators, coming from the two types of rigid configurations of 6 points in $\RP^2$.
\end{Th}
\begin{pf} As before, it is enough to consider 6-tuples of line segments $L=(L_1,\dots, L_6)$, where
$L_i=[0,u_i]$ for some $u_i\in\R^3$, $i\in[6]$. In this case $\V_I(L)=|\det(u_I)|$, where $u_I$ is the matrix whose columns are $u_i$ for $i\in I$.

Let $x_I=\det(u_I)\det(u_{I^c})$. Fix any $J\subset[6]$ of size two and $\ell\in J^c$. Without loss of generality we may assume $J=1\,2$ and $\ell=6$.
Then the inequalities \re{6-zonoids-in-3d} are
\begin{equation}\label{e:6-zonoids-in-3d-explicit}
%0\leq \sum_{k\not\in I\supset J}\!\!\vv_{I}-\vv_{J\cup \{k}\leq \sum_{I\not > J}\vv_I,
\vv_{1\,2\,6}\leq \vv_{1\,2\,3}+\vv_{1\,2\,4}+\vv_{1\,2\,5}\leq \vv_{1\,2\,6} +\sum_{1\,2\,\not\subset I\ni 6}\!\vv_I.
\end{equation}
According to the Grassmann-Plücker relation, we have 
\begin{equation}\label{e:Plucker}
x_{1\,2\,3}-x_{1\,2\,4}+x_{1\,2\,5}-x_{1\,2\,6}=0.
\end{equation}
Since $\vv_I=|x_I|$, the triangle inequality provides
\begin{equation}\label{e:abs-Plucker}
%\vv_{1\,2\,6}=|x_{1\,2\,6}|=|x_{1\,2\,3}-x_{1\,2\,4}+x_{1\,2\,5}|\leq |x_{1\,2\,3}|+|x_{1\,2\,4}|+|x_{1\,2\,5}|=\vv_{1\,2\,3}+\vv_{1\,2\,4}+\vv_{1\,2\,5}.\nonumber
\vv_{1\,2\,6}=|x_{1\,2\,3}-x_{1\,2\,4}+x_{1\,2\,5}|\leq \vv_{1\,2\,3}+\vv_{1\,2\,4}+\vv_{1\,2\,5}.\nonumber
\end{equation}
This proves the first inequality in \re{6-zonoids-in-3d-explicit}. 

Now we prove the second inequality in  \re{6-zonoids-in-3d-explicit}. %we need to show
%\begin{equation}\label{e:2}
%\vv_{1\,2\,3}+\vv_{1\,2\,4}+\vv_{1\,2\,5}\leq \vv_{1\,2\,6} +\sum_{1\,2\,\not\subset I\ni 6}\!\vv_I
%\end{equation}
For this we claim that we can express every sum $\pm x_{1\,2\,3}\pm x_{1\,2\,4}\pm x_{1\,2\,5}$
as a linear combination of $x_{1\,2\,6}$ and $\left\{x_I :  1\,2\,\not\subset I\ni 6\right\}$ with coefficients in $\{-1,0,1\}$. Then, applying the triangle inequality as above, we obtain the second inequality in \re{6-zonoids-in-3d-explicit}:
$$
\vv_{1\,2\,3}+\vv_{1\,2\,4}+\vv_{1\,2\,5}=\max\left\{\pm x_{1\,2\,3}\pm x_{1\,2\,4}\pm x_{1\,2\,5}\right\}\leq
\vv_{1\,2\,6}+\sum_{1\,2\,\not\subset I\ni 6}\!\vv_I.$$
To prove the claim, first note that one of such sums is $x_{1\,2\,3}- x_{1\,2\,4}+ x_{1\,2\,5}$ which is simply  $x_{1\,2\,6}$
by \re{Plucker}. Up to sign, there are three more sums to consider
%Thus, to prove the claim it is enough to consider the sums
$$x_{1\,2\,3}+ x_{1\,2\,4}+ x_{1\,2\,5},\quad x_{1\,2\,3}+ x_{1\,2\,4}- x_{1\,2\,5},
\ \text{ and }\ x_{1\,2\,3}- x_{1\,2\,4}- x_{1\,2\,5}.$$ 
We will express the first sum; the other two cases are similar. We have the following Grassmann-Plücker relations
\begin{align}%\label{e:Pluckers}
x_{1\,2\,4}-x_{1\,3\,4}-x_{1\,4\,5}+x_{1\,4\,6}=0,\nonumber\\
x_{1\,2\,4}+x_{2\,3\,4}+x_{2\,4\,5}-x_{2\,4\,6}=0.\nonumber
\end{align}
Adding them together with \re{Plucker} produces
$$x_{1\,2\,3}+x_{1\,2\,4}+x_{1\,2\,5}=x_{1\,2\,6}+x_{1\,3\,4}+x_{1\,4\,5}-x_{1\,4\,6}-x_{2\,3\,4}-x_{2\,4\,5}+x_{2\,4\,6}.$$
Since $x_I=x_{I^c}$ we can replace appropriate index subsets with their complements, hence, 
$$x_{1\,2\,3}+x_{1\,2\,4}+x_{1\,2\,5}=x_{1\,2\,6}+x_{2\,5\,6}+x_{2\,3\,6}-x_{1\,4\,6}-x_{1\,5\,6}-x_{1\,3\,6}+x_{2\,4\,6},$$
as required.
%Let $x_I=\det(u_I)\det(u_{I^c})$. Fix any $J\subset[6]$ of size two and let $i_1<i_2<i_3<i_4$ be the four numbers in the complement $J^c$. Then, according to the Grassmann-Plücker relation, we have 
%\begin{equation}\label{e:Plucker}
%\sum_{\ell=1}^4(-1)^\ell x_{J\cup \{i_\ell}=0,\nonumber
%\end{equation}
\end{pf}

\section{Non-polyhedrality result for $C_{n,2}$} 
\label{S:nonpolyhedrality} 

In general, the polyhedrality of $C_{n,d}$ is not expected. 
The smallest case when the cone is not polyhedral is $(n,d)=(8,2)$, as we show in the next theorem. Our proof uses some computations done in SageMath, see
\url{https://github.com/jsimonrichard/mixed-volumes-n-6-d-2-3}  for the detailed code.

\begin{Th}\label{T:nonpolyhedral}
	The cone $C_{2k,2}$ is not polyhedral for each integer $k \ge 4$. 
\end{Th}

\begin{pf}
	It suffices to verify non-polyhedrality for $C_{8,2}$, as $C_{2k,2}$ with $k>4$ can be mapped onto $C_{8,2}$ via a coordinate projection that keeps the coordinates $\V_{I_1 | \cdots | I_k}$ with the sets $I_5,\ldots, I_k$ being fixed. For example, to project $C_{10,2}$ onto $C_{8,2}$, we may keep all $V_{ I_1 | I_2 | I_3 | I_4 |  \{ 9 , 10 \} }$, where $\{9,10\}$ is the fixed part of the $2+2+2+2+2$ partition of $[10]$ and $I_1 | \cdots | I_4$ runs over $\Prt(2,8)$. Thus, we deal with the cone $C_{8,2} \subset \R^{105}$. 
	
	The classification of rigid configuration of eight points in $\RP^1$ up to projective transformations is as follows: (a) a pair of points with multiplicity four $4\cdot\{0,1\}$, (b)
	a pair of double points and a multiplicity four point $\{2\cdot 0, 2\cdot 1, 4\cdot\infty\}$,
	(c) a pair of triple points and a double point $\{3\cdot 0, 3\cdot 1, 2\cdot\infty\}$, and (d)
	a 1-parameter family of four double points $2\cdot \{0, 1, t, \infty\}$, for $t\in\RP^1$. Configurations with points of multiplicity larger than four produce the zero vector under the map $\Phi$, hence, we ignore them.
	
	First, we compute the cone $C_0\subset \R^{105}$ generated by $\Phi(U)$ over all 8-tuples $U$ corresponding to finite configurations (a)--(c). This cone is polyhedral of dimension 91, generated by 315 rays. For an 8-tuple $U$ corresponding to $2\cdot \{0, 1, t, \infty\}$, the image
	$\Phi(U)$ is a piecewise quadratic curve. More explicitly, on each of the three intervals $[0,1]$, $[1,\infty]$, and $[-\infty, 0]$, we have $\Phi(U)=c_0+c_1t+c_2t^2$ for some linearly independent vectors $c_0,c_1,c_2\in\R^{105}$. Taking all possible $U$
	corresponding to $2\cdot \{0, 1, t, \infty\}$ produces 630 such curves. In fact, some of the permutations
	amount to permuting the three pieces of each of the piecewise quadratic curve. It follows that one
	can restrict the computation to only pieces defined on one of the segments, say, $[0,1]$.
	Therefore, the cone $C_{8,2}$ is the conic hull of $C_0$ together with the union of 630 
	quadratic curves defined on $[0,1]$.
	
	As $C_{8,2}$ is contained in the positive orthant, we can ``dehomogenize'' it by 
	considering the affine slice $K=C_{8,2}\cap\{x\in\R^{105} : \sum_i x_i=1\}$.
	Note that $K$ is the convex hull of normalized points $\Phi(U)/\sigma(\Phi(U))$, where we denote
	$\sigma(x)=\sum_{i} x_i$. Thus $C_{8,2}$ is polyhedral if and only if $K$ is a polytope.
	
	Next project $K$ onto a 2-dimensional plane using $\pi:\R^{105}\to\R^2$, 
	$\pi(v)=(\overline{m}\cdot v, \overline{n}\cdot v)$ for some $\overline{m},\bar{n}\in\R^{105}$. In our implementation we chose
	$\overline{m}$ generic and $\bar{n}$ orthogonal to $\spn(C_0)$.  Here is our particular choice of the coordinates of $\overline{m}$ and $\overline{n}$: 
	{ \scriptsize 
		\begin{align*}
		\begin{array}{l|rrrrrrrrrrrrrrrrrrrr}
		i & 1 & 2 & 3 & 4 & 5 & 6 & 7 & 8 & 9 & 10 & 11 & 12 & 13 & 14 & 15 & 16 & 17 & 18 & 19 & 20 
		\\ \hline \hline 
		\overline{m}_i & 5 & 5 & 4 & 5 & 3 & 2 & 5 & 4 & 0 & 4 & 2 & 4 & 5 & 3 & 1 & 2 & 5 & 2 & 5 & 4 \\
		\overline{n}_i & 0 & 0 & 0 & 0 & 0 & 0 & 0 & 0 & 1 & -1 & 1 & -1 & 0 & 0 & 1 & -1 & 1 & -1 & 0 & 0
		\\ \hline 
		i & 21 & 22 & 23 & 24 & 25 & 26 & 27 & 28 & 29 & 30 & 31 & 32 & 33 & 34 & 35 & 36 & 37 & 38 & 39 & 40 
		\\ \hline \hline 
		\overline{m}_i & 0 & 4 & 5 & 4 & 6 & 3 & 5 & 3 & 0 & 4 & 1 & 2 & 5 & 3 & 1 & 6 & 1 & 6 & 4 & 3 \\
		\overline{n}_i & 1 & -1 & 1 & -1 & 0 & 0 & 0 & 0 & 0 & 0 & 0 & 0 & 0 & 0 & 0 & 0 & 0 & 0 & 0 & 0
		\\ \hline 
		i & 41 & 42 & 43 & 44 & 45 & 46 & 47 & 48 & 49 & 50 & 51 & 52 & 53 & 54 & 55 & 56 & 57 & 58 & 59 & 60 
		\\ \hline \hline 
		\overline{m}_i & 4 & 6 & 0 & 5 & 4 & 4 & 4 & 5 & 5 & 1 & 4 & 0 & 1 & 3 & 2 & 4 & 5 & 1 & 1 & 4 \\
		\overline{n}_i & 0 & 0 & 0 & 0 & 0 & 0 & 1 & -1 & 1 & -1 & 0 & 0 & 1 & -1 & 1 & -1 & 0 & 0 & 1 & -1 \\ 
		\hline 
		i & 61 & 62 & 63 & 64 & 65 & 66 & 67 & 68 & 69 & 70 & 71 & 72 & 73 & 74 & 75 & 76 & 77 & 78 & 79 & 80 
		\\ \hline \hline 
		\overline{m}_i & 
		0 & 5 & 4 & 1 & 0 & 4 & 6 & 3 & 5 & 5 & 5 & 3 & 1 & 6 & 2 & 3 & 1 & 0 & 3 & 2 \\
		\overline{n}_i & 1 & -1 & 0 & 0 & 0 & 0 & 0 & 0 & 0 & 0 & 0 & 0 & 0 & 0 & 0 & 0 & 0 & 0 & 0 & 0 \\  \hline 
		i & 81 & 82 & 83 & 84 & 85 & 86 & 87 & 88 & 89 & 90 & 91 & 92 & 93 & 94 & 95 & 96 & 97 & 98 & 99 & 100 \\
		\hline \hline 
		\overline{m}_i & 5 & 3 & 3 & 4 & 6 & 6 & 4 & 4 & 1 & 0 & 4 & 6 & 3 & 3 & 5 & 4 & 0 & 2 & 4 & 5 \\
		\overline{n}_i & 
		0 & 0 & 0 & 0 & 0 & 0 & 0 & 0 & 0 & 0 & 0 & 0 & 0 & 0 & 0 & 0 & 0 & 0 & 0 & 0 \\ \hline 
		i & 101 & 102 & 103 & 104 & 105 \\
		\hline \hline 
		\overline{m}_i & 5 & 3 & 8 & 0 & 4 \\
		\overline{n}_i & 0 & 0 & 0 & 0 & 0
		\end{array}
		\end{align*}
	}
	
	We show that the projected 2-dimensional body $\pi(K)$ is not a polygon and, hence, $K$
	is not a polytope.
	
	We compute $\pi(K)$ by first projecting the normalized generators and normalized curves 
	and then taking the convex hull.
	Since $n$ is orthogonal to $\spn(C_0)$, the projection of the normalized generators
	of $C_0$ lie on the $x$-axis in $\R^2$, see \rf{many_curves}. The two extremal
	points in the projection are $(11/4,0)$ and $(143/36,0)$.
	
	\begin{figure}[h]
		\begin{center}
			\includegraphics[scale=.5]{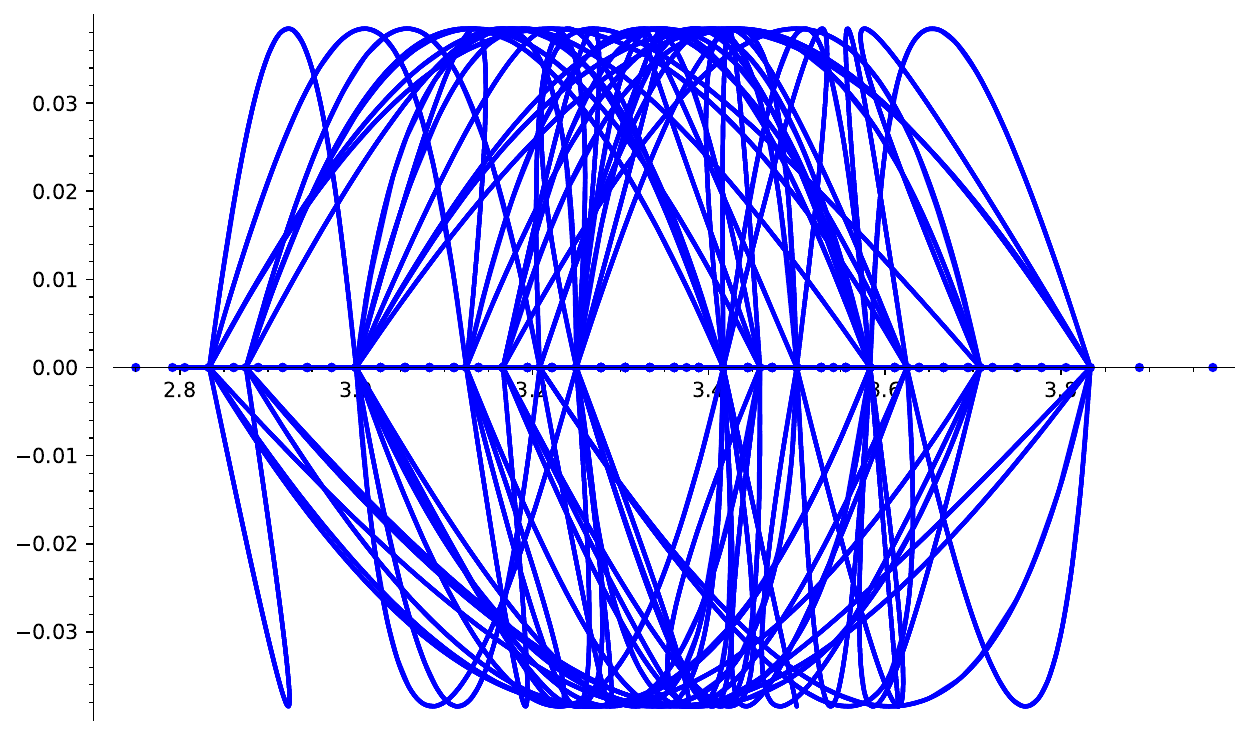}
		\end{center}
		\caption{The projection of normalized generators and normalized curves}
		\label{F:many_curves}
	\end{figure}
	
	The 630 normalized quadratic curves are rational curves. Their projections are depicted in \rf{many_curves} in blue. Note that some of the curves lie in $\spn(C_0)$ and, hence, also project onto the $x$-axis.
	There are four extremal rational curves, depicted in Fig.~\ref{F:convex_hull}, in the projection $\{(x_i(t),y_i(t)) : t\in[0,1], 1\leq i\leq 4\}$, with the parametrizations given as follows: 
	\begin{align*}
	x_1(t) &= \frac{6t^2 - 7t - 17}{2(t^2 - t - 3)},\quad
	x_2(t) = \frac{30t^2 - 31t - 68}{8(t^2 - t - 3)},\quad
	x_3(t) = \frac{12t^2 - 9t - 46}{4(t^2 - t - 3)},\\
	x_4(t) &= \frac{33t^2 - 28t - 92}{8(t^2 - t - 3)},\quad\quad y_i(t) = (-1)^{i-1}\frac{t^2-t}{2(t^2 - t - 3)}\
	1\leq i\leq 4.
	\end{align*}
	For each of the four displayed rational curves, we compute the tangent lines from one of the two extremal points $(11/4,0)$ and $(143/36,0)$. The corresponding tangency points have the following $y$-coordinates:
	$$
	\frac{7}{2}-2\sqrt{3},\quad
	-\frac{7}{23}+\frac{5\sqrt{6}}{46},\quad
	\frac{245}{181}-\frac{57\sqrt{70}}{362},\quad
	-\frac{5}{7}+\frac{3\sqrt{10}}{14}.
	$$
	A direct calculation shows that the first and third values lie strictly between $0$ and the maximum value $1/26$, while the second and fourth lie strictly between the minimum value $-1/26$ and $0$. %Thus the tangencies occur at interior points of the four curves. 
	Since the tangent lines are supporting lines for $\pi(K)$, there are nontrivial arcs of these rational curves that lie on the boundary of $\pi(K)$. Therefore, $\pi(K)$ is not a polygon. We depict the extremal curves together with the tangency points and the boundary of $\pi(K)$ in \rf{convex_hull}. 
	
	\begin{figure}[h]
		\begin{center}
			\includegraphics[scale=.5]{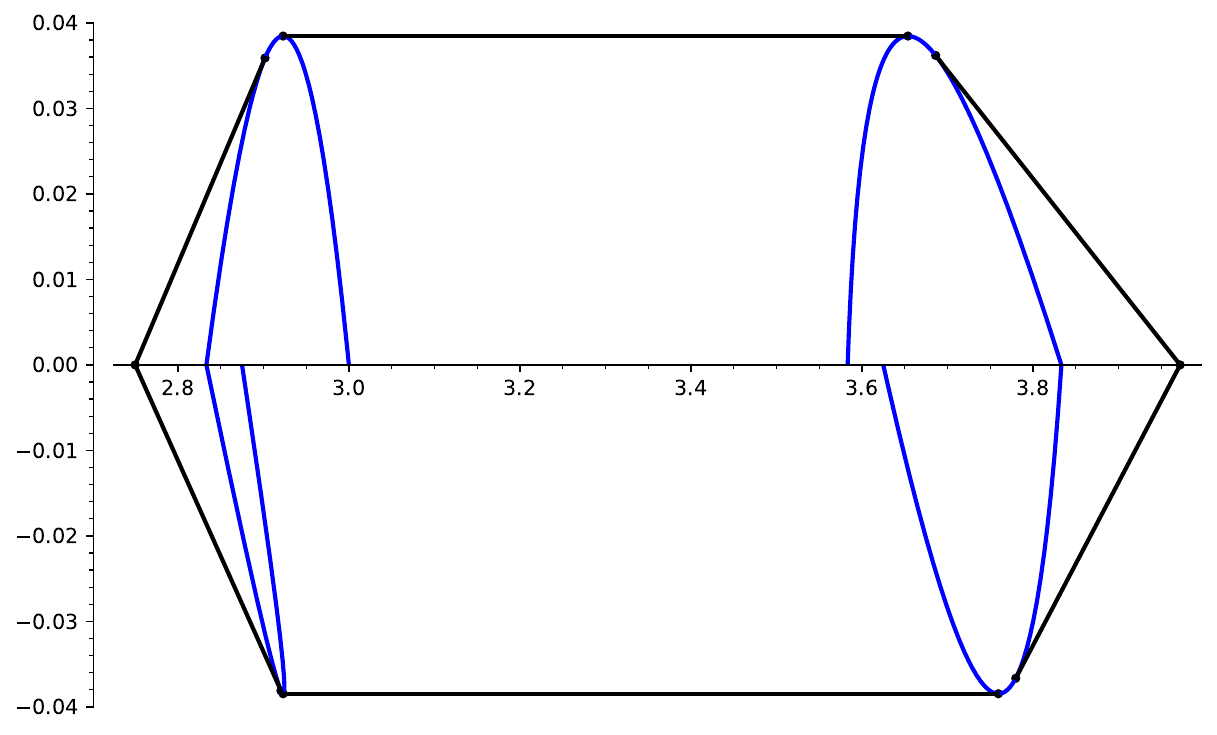}
		\end{center}
		\caption{The projection $\pi(K)$}
		\label{F:convex_hull}
	\end{figure}
	
\end{pf}

\begin{Rem} 
	The results of this section show that,  in dimension two, within the family of centrally-symmetric convex bodies, we have a transition from polyhedrality to non-polyhedrality of $C_{n,2}$, when $n$ grows. 
	We expect that such a transition also occurs for the cones $C_{n,d}$ in every dimension $d \ge 3$. Studying this phenomenon is a subject of future research. 
\end{Rem}

\subsection*{Data availability statement} 
No data sets were generated or analyzed during the current study.

\bibliographystyle{plain} 
\bibliography{lit}

\begin{thebibliography}{10}

\bibitem{adiprasito2017hodge}
Karim Adiprasito, June Huh, and Eric Katz.
\newblock Hodge theory of matroids.
\newblock {\em Notices of the AMS}, 64(1), 2017.

\bibitem{adiprasito2018hodge}
Karim Adiprasito, June Huh, and Eric Katz.
\newblock Hodge theory for combinatorial geometries.
\newblock {\em Annals of Mathematics}, 188(2):381--452, 2018.

\bibitem{AS23}
Gennadiy Averkov and Ivan Soprunov.
\newblock Pl{\"u}cker-type inequalities for mixed areas and intersection
  numbers of curve arrangements.
\newblock {\em International Mathematics Research Notices},
  2023(18):16015--16050, 2023.

\bibitem{bochnak2013real}
Jacek Bochnak, Michel Coste, and Marie-Fran{\c{c}}oise Roy.
\newblock {\em Real algebraic geometry}, volume~36.
\newblock Springer Science \& Business Media, 2013.

\bibitem{branden2020lorentzian}
Petter Br{\"a}nd{\'e}n and June Huh.
\newblock Lorentzian polynomials.
\newblock {\em Annals of Mathematics}, 192(3):821--891, 2020.

\bibitem{ZonoidAlgebra}
Paul Breiding, Peter Bürgisser, Antonio Lerario, and Léo Mathis.
\newblock The zonoid algebra, generalized mixed volumes, and random
  determinants.
\newblock {\em Advances in Mathematics}, 402:108361, 2022.

\bibitem{FZ23}
Matthieu Fradelizi, Mokshay Madiman, Mathieu Meyer, and Artem Zvavitch.
\newblock On the volume of the {M}inkowski sum of zonoids.
\newblock {\em J. Funct. Anal.}, 286(3):Paper No. 110247, 41, 2024.

\bibitem{Heine38}
Rudolf Heine.
\newblock Der {W}ertvorrat der gemischten {I}nhalte von zwei, drei und vier
  ebenen {E}ibereichen.
\newblock {\em Mathematische Annalen}, 115(1):115--129, 1938.

\bibitem{huang2025realizations}
Daoji Huang, June Huh, Mateusz Micha{\l}ek, Botong Wang, and Shouda Wang.
\newblock Realizations of homology classes and projection areas.
\newblock {\em arXiv preprint arXiv:2505.08881}, 2025.

\bibitem{khovanskii1978newton}
Askold~G Khovanskii.
\newblock Newton polyhedra and the genus of complete intersections.
\newblock {\em Functional Analysis and its applications}, 12(1):38--46, 1978.

\bibitem{lovasz1987matching}
L\'aszl\'o Lov\'asz.
\newblock Matching structure and the matching lattice.
\newblock {\em J. Combin. Theory Ser. B}, 43(2):187--222, 1987.

\bibitem{postnikov2018positive}
Alexander Postnikov.
\newblock Positive {G}rassmannian and polyhedral subdivisions.
\newblock In {\em Proceedings of the International Congress of Mathematicians:
  Rio de Janeiro 2018}, pages 3181--3211. World Scientific, 2018.

\bibitem{scheiderer2024course}
Claus Scheiderer.
\newblock {\em A course in real algebraic geometry}.
\newblock Springer, 2024.

\bibitem{Schneider2014}
Rolf Schneider.
\newblock {\em Convex bodies: the {B}runn-{M}inkowski theory}, volume 151 of
  {\em Encyclopedia of Mathematics and its Applications}.
\newblock Cambridge University Press, Cambridge, expanded edition, 2014.

\bibitem{Shephard60}
GC~Shephard.
\newblock Inequalities between mixed volumes of convex sets.
\newblock {\em Mathematika}, 7(2):125--138, 1960.

\bibitem{speyer2004tropical}
David Speyer and Bernd Sturmfels.
\newblock The tropical {G}rassmannian.
\newblock {\em Adv. Geom.}, 4(3):389--411, 2004.

\bibitem{Stan}
Richard~P Stanley.
\newblock Two combinatorial applications of the {A}leksandrov-{F}enchel
  inequalities.
\newblock {\em Journal of Combinatorial Theory, Series A}, 31(1):56--65, 1981.

\bibitem{teissier1979theoreme}
Bernard Teissier.
\newblock Du th{\'e}oreme de l’index de hodge aux in{\'e}galit{\'e}s
  isop{\'e}rim{\'e}triques.
\newblock {\em CR Acad. Sci. Paris S{\'e}r. AB}, 288(4):A287--A289, 1979.

\bibitem{sagemath}
{The Sage Developers}.
\newblock {\em {S}ageMath, the {S}age {M}athematics {S}oftware {S}ystem
  ({V}ersion x.y.z)}, 2024.
\newblock {\tt https://www.sagemath.org}.

\bibitem{williams2021positive}
Lauren~K. Williams.
\newblock The positive {G}rassmannian, the amplituhedron, and cluster algebras.
\newblock In {\em I{CM}---{I}nternational {C}ongress of {M}athematicians.
  {V}ol. 6. {S}ections 12--14}, pages 4710--4737. EMS Press, Berlin, [2023]
  \copyright 2023.

\end{thebibliography}

\end{document}